\documentclass{m2an}
\usepackage{bm}
\usepackage{amsmath}
\usepackage{amssymb}
\usepackage{amsthm}
\usepackage{enumerate}
\usepackage{float}
\usepackage{tikz}
\usepackage{epstopdf}
\usepackage{booktabs}
\usepackage{makecell}
\usepackage{hyperref}
\begin{document}
\title{Computing Both Upper and Lower Eigenvalue Bounds by HDG Methods}
\author{Qigang Liang}\address{School of Mathematical Science, Tongji University, Shanghai 200092, China;\\
\email{qigang$\_$liang@tongji.edu.cn; xuxj@tongji.edu.cn; lyyuan@tongji.edu.cn}}
\author{Xuejun Xu}\sameaddress{1}
\author{Liuyao Yuan}\sameaddress{1}
%
%
\begin{abstract}
In this paper, we observe an interesting phenomenon for a hybridizable discontinuous Galerkin (HDG) method for eigenvalue problems. Specifically, using the same finite element method, we may achieve both upper and lower eigenvalue bounds simultaneously, simply by the fine tuning of the stabilization parameter. Based on this observation, a high accuracy algorithm for computing eigenvalues is designed to yield higher convergence rate at a lower computational cost. Meanwhile, we demonstrate that certain type of HDG methods can only provide upper bounds. As a by-product, the asymptotic upper bound property of the Brezzi-Douglas-Marini mixed finite element is also established. Numerical results supporting our theory are given.
\end{abstract}
%
%
\subjclass{65N12, 65N15, 65N30}
\keywords{PDE eigenvalue problems, HDG methods, upper and lower bound, high accuracy}
\maketitle
\section*{Introduction}

Eigenvalue problems play an important role in various fields such as quantum mechanics, fluid mechanics, stochastic processes, etc. (see, e.g., \cite{MR2787541,MR3617375,MR3499456}). Since it is  almost impossible to obtain the exact eigenvalues, finding efficient numerical approximations is a fundamental and long-standing work. It is well known that conforming finite element methods provide guaranteed upper bounds based on the Courant-Fischer principle. Nonconforming finite element methods and mixed finite element methods are alternative possible ways of providing lower bounds. In this paper, we investigate the lower and upper bound properties of mixed Laplacian eigenvalue problems discretized by the hybridizable discontinuous Galerkin (HDG) methods.

There are many results to study the lower bound property of discrete eigenvalues by nonconforming elements. Boffi \cite{MR2652780} presents several examples on numerical computations of eigenvalue problems using nonconforming and mixed finite element methods. Armentano and Dur\'{a}n \cite{MR2040799} make an important advance regarding lower bounds for eigenvalue problems by the Crouzeix-Raviart (CR) element. It is proved that the discrete eigenvalues provided by the CR element are smaller than exact ones when the corresponding eigenfunctions belong to $H^{1+r}(\Omega) \cap H_0^1(\Omega)$ with $0<r<1$. Hu et al. \cite{MR3254372} introduce a new systematic method for constructing nonconforming finite elements to obtain lower bounds of exact eigenvalues, including the enriched nonconforming rotated $Q_1$ element and the Wilson element. Furthermore, conforming approximations of exact eigenpairs are constructed (see \cite{MR3385147}) based on the nonconforming discrete eigenfunctions. Carstensen et al. \cite{MR4048619} establish a guaranteed lower bound (GLB) criterion for eigenvalue approximations by a HDG method with Lehrenfeld–Sch{\"o}berl stabilization. This approach was later generalized to the framework of hybrid high-order (HHO) methods in \cite{MR4332791} and the bi-Laplacian eigenvalue problem in \cite{MR4592002}. Lin et al. \cite{MR2429874} apply the integral-identity technique to obtain the asymptotic expansion of discrete eigenvalues by nonconforming finite element methods, and lower bound properties are achieved as a by-product (also see \cite{MR2427427,MR2490184}).

However, there has been limited research on analyzing the upper or lower bound properties of eigenvalues using mixed finite element methods. Some numerical examples demonstrate that conforming mixed finite elements do not necessarily lead to upper bounds for eigenvalues. For instance, numerical results in \cite{MR2400623} show that the mixed Bernadi-Raugel element for the Stokes eigenvalue problem produce lower bounds of exact eigenvalues. Additionally, numerical results in \cite{MR2652780} demonstrate that, for the mixed Laplacian eigenvalue problem, the eigenvalues computed by the Raviart-Thomas (RT) element may approximate the exact ones from below or above on some unstructured meshes, while all eigenvalues are approximated from above on specific structured meshes of squares. To the best of our knowledge, there is no theoretical proof for these phenomena yet. For the theory of lower or upper bounds of eigenvalues by the mixed finite element method, one early result is achieved in \cite{MR2768711} for the Stokes eigenvalue problem. The paper introduces a suitable theoretical framework to analyze the lower bound property of the eigenvalues of the Stokes operator with nonconforming mixed finite elements, including the (enriched) CR element and the (enriched) rotated $Q_1$ element. Later, Yang et al. \cite{MR3063973} propose several mixed finite element identities to analyze the upper bound property of discrete eigenvalues using the MINI element, etc. Most recently, Gallistl \cite{MR4570331} employs mixed finite element methods to provide computable GLBs for the Laplacian, linear elasticity, and the Steklov eigenvalue problem through a post-processing of the discrete eigenvalues.

In recent decades, the HDG methods have gained popularity in solving partial differential equations due to their three well-known advantages: reduced system size, flexible stabilization, and enhanced accuracy through post-processing. Firstly, by introducing an additional unknown variable $\lambda_h$ (interpreted as a Lagrange multiplier), a ``static condensation'' technique \cite{MR1115205} can be applied to eliminate internal degrees of freedom, resulting in a smaller stiffness matrix compared to the one associated to the original variables $(\mathbf p_h,u_h)$. This advantage becomes particularly significant when dealing with high order elements. Secondly, the stabilization of HDG methods with a penalty term provides a highly flexible way to select suitable approximate spaces. Thirdly, a post-processing technique can be applied to obtain higher order solutions (see \cite{MR2629996,MR3315499}). The idea of hybrid methods can be tracked back to the 1960s \cite{MR3618552}. A unifying framework for the hybridization of discontinuous Galerkin methods is introduced in \cite{MR2485455} by Cockburn et al. A new  technique for error analysis of HDG methods is put forward in \cite{MR2629996}. Recently, Hong et al. \cite{MR4191130} provide a number of new estimates on the stability and convergence of both HDG and weak Galerkin (WG) methods based on the Babu\v{s}ka-Brezzi theory. For eigenvalue problems, hybridization and post-processing techniques for the RT approximation of second order elliptic eigenvalue problems are presented in \cite{MR2669393} (known as the HRT method). The approximation by the HDG method for the Laplacian eigenvalue problem is considered in \cite{MR3315499}, where the convergence rates of approximate eigenvalues and the corresponding eigenfunctions are obtained. It should be noted that the aforementioned advantages of HDG finite element discretization for source problems can also be inherited by eigenvalue problems.

In this paper, we study two types of HDG approximations of discrete elliptic eigenvalue problems, and present two different approximation results for the upper and lower bounds of eigenvalues. Specifically, the particular HDG methods we consider are referred to as the LDG-H method in \cite{MR2485455}.  Recently, we have also used a special WG method (distinct from HDG methods) to do similar works \cite{MR4476865}, which requires continuity constraints along the interelement boundaries. In comparison, the work of HDG here is more general. We must emphasize that the work of HDG is nontrivial. Firstly, the HDG methods no longer restrict continuity along interelement boundaries. Secondly, to obtain convergence results for eigenvalue problems, we have to introduce some new estimates on the convergence of the corresponding source problems discretized by the HDG methods.

For the first type of HDG methods, it is proved that both lower and upper bounds of eigenvalues could be obtained merely by adjusting the penalty parameter in the bilinear form, that is, if the penalty parameter is sufficiently small (resp. large), the discrete eigenvalue is a lower (resp. upper) bound of the corresponding exact eigenvalue. We propose two identities to analyze the upper and lower bound properties. A post-processing technique for computing eigenvalues is then applied to obtain high accuracy approximations.  For the second type of HDG methods, we prove that the discrete eigenvalues approximate the exact ones from above, regardless of how the penalty parameter changes. This is achieved by analyzing the superconvergence result between the discrete eigenfunction and the projection which corresponds to the exact eigenfunction. In addition, the upper bound property of eigenvalues computed by the Brezzi-Douglas-Marini (BDM) element is achieved as a by-product. To the best knowledge of the authors, this result never appears in the literature.

The outline of this paper is organized as follows: In Section \ref{sec:Model Problem and Preliminaries},  the model problem and some preliminaries are introduced. In Section \ref{sec:Gradient-based HDG discretization}, we explore upper and lower bound properties of eigenvalues by the first type of HDG methods (gradient-based). The second type of HDG methods (divergence-based), which admits upper bounds of eigenvalues, is analyzed in Section \ref{sec:Divergence-based HDG discretization}. A post-processing technique is designed in Section \ref{sec:A High Accuracy Algorithm}. In Section \ref{sec:Numerical Experiments}, we present some numerical examples to support theoretical findings. Finally, some conclusions are given in Section \ref{sec:Conclusions}.
\section{Model Problem and Preliminaries}\label{sec:Model Problem and Preliminaries}
In this section, we first describe some basic notations in Subsection \ref{subs:Notations} and introduce the Laplacian eigenvalue problem in Subsection \ref{subs:Model Problem}.  Then the corresponding HDG finite element discretization scheme is presented in Subsection \ref{subs:HDG finite element discretization scheme}.

\subsection{Notations}
\label{subs:Notations}
Throughout this paper, we use standard notations for Sobolev spaces $H^{m}(D)$ and $H_{0}^{m}(D)$ with their associated norms $\|\cdot\|_{m,D}$ and semi-norms $|\cdot|_{m,D}$ (see \cite{MR0450957}). If $D = \Omega$, we use the abbreviated notations $\|\cdot\|_{m}$ and $|\cdot|_{m}$ for simplicity. We denote by $L^{2}(D):=H^{0}(D)$, and the $L^2$-inner products on $D$ and $\partial D$ are denoted by $(\cdot,\cdot)_{D}$ and $\left\langle \cdot,\cdot\right\rangle_{\partial D}$, respectively. We shall drop the subscript $D$ in the product notation $(\cdot,\cdot)_{D}$ if $D = \Omega$. Moreover, the norms of $L^2(D)$ and $L^2(\partial D)$ are denoted by $\|\cdot\|_{0,D}$ and $\|\cdot\|_{0,\partial D}$, respectively. For vector field space, we use bold font $\mathbf{L}^{2}(\Omega)$ and $\mathbf{H}^{m}(\Omega)$ to represent $[L^{2}(\Omega)]^{2}$ and $[H^{m}(\Omega)]^{2}$, respectively. We also denote by 
$$\mathbf H(\operatorname{div}, \Omega) := \{\mathbf q \in \mathbf L^2(\Omega) : \operatorname{div} \mathbf q \in L^2(\Omega)\}$$ equipped with the norm  $\|\mathbf q\|_{\operatorname{div}}^2=(\mathbf q, \mathbf q)+(\operatorname{div} \mathbf q, \operatorname{div} \mathbf q)$. We use letter $C$, with or without subscripts, to denote a generic constant, independent of the mesh size, which stands for different values at different occurrences. Additionally, we note that the dependencies on the penalty parameter will always be explicitly mentioned.

Let $\Omega$ be a convex polygonal domain in $\mathbb{R}^2$ with boundary $\partial \Omega$ and $\mathcal{T}_{h}$ be a shape-regular and quasi-uniform triangulation of $\Omega$. Define $h := \max_{K \in \mathcal{T}_{h}} h_K$, where $h_K$ denotes the diameter of $K$, $K \in \mathcal{T}_{h}$.  We denote $\mathcal{E}_h^i$ the set of interior edges in $\mathcal{T}_{h}$, and $\mathcal{E}_h^\partial$ the set of boundary edges. Furthermore, we denote by $\mathcal{E}_h := \mathcal{E}_h^i \cup \mathcal{E}_h^\partial$ the set of all edges in the triangulation. For any edge $e \in \mathcal{E}_h$, $h_e$ denotes its diameter. For an interior edge $e \in \mathcal{E}_h^i$ that is the common edge of two adjacent triangles $K^\pm$, and for any vector-valued function $\mathbf q$, we define the jump on $e$ as follows:
\begin{equation*}
[\mathbf q]:=\left.\mathbf q\right|_{K^+} \cdot \mathbf n^{+}+\left.\mathbf q\right|_{K^-} \cdot \mathbf n^{-},
\end{equation*}
where $\mathbf n^+$ and $\mathbf n^-$ are the outward unit normal vectors corresponding to $\partial K^+$ and $\partial K^-$, respectively. If $e \subset \partial \Omega$, we denote by $\mathbf n$ the unit normal vector pointing outside $\Omega$ and define $[\mathbf q]:=\mathbf q \cdot \mathbf n$.

Finally, we introduce some inner products and norms as follows:
\begin{alignat*}{3}
(\cdot, \cdot)_{\mathcal{T}_h} &:= \sum_{K \in \mathcal{T}_h}(\cdot, \cdot)_K,&\qquad \langle\cdot, \cdot\rangle_{\partial \mathcal{T}_h} &:= \sum_{K \in \mathcal{T}_h}\langle\cdot, \cdot\rangle_{\partial K}, \\
\langle\cdot, \cdot\rangle_{\mathcal{E}_h} &:= \sum_{e \in \mathcal{E}_h}\langle\cdot, \cdot\rangle_e, & \langle\cdot,  \cdot\rangle_{\mathcal{E}_h^i} &:= \sum_{e \in \mathcal{E}_h^i}\langle\cdot, \cdot\rangle_e, \\
\|\cdot\|_{\partial \mathcal{T}_h}^2 &:= \langle\cdot, \cdot\rangle_{\partial \mathcal{T}_h},& \|\cdot\|_{\mathcal{E}_h^i}^2 &:= \langle\cdot, \cdot\rangle_{\mathcal{E}_h^i}.
\end{alignat*}
We now give more details for $\langle \cdot, \cdot\rangle_{\partial \mathcal{T}_h}$. For any scalar-valued function $v$ and vector-valued function $\mathbf q$,
\begin{equation*}
\langle v, \mathbf q \cdot \mathbf n\rangle_{\partial \mathcal{T}_h} 
= \langle v, \mathbf q \cdot \mathbf n_K\rangle_{\partial \mathcal{T}_h}
=\sum_{K \in \mathcal{T}_h}\langle v, \mathbf q \cdot \mathbf n\rangle_{\partial K}
=\sum_{K \in \mathcal{T}_h}\left\langle v, \mathbf q \cdot \mathbf n_K\right\rangle_{\partial K},
\end{equation*}
where $\mathbf n_K$ is the outward unit normal vector along $\partial K$.

\subsection{Model Problem}
\label{subs:Model Problem}
Consider the Laplacian eigenvalue problem and set $\mathbf p=-\nabla u$ to obtain the following mixed form:
\begin{equation}\label{eq:EVP}
\begin{cases}
\mathbf p+\nabla u=0 & \text { in } \Omega, \\
\operatorname{div} \mathbf p=\lambda u & \text { in } \Omega, \\
u = 0 & \text { on } \partial\Omega.
\end{cases}
\end{equation}
Let $\mathbf Q := \mathbf H(\operatorname{div},\Omega)$ and $V := L^2(\Omega)$. The weak form of \eqref{eq:EVP} reads: Find $(\lambda, \mathbf p, u) \in \mathbb{R} \times \mathbf Q \times V$, such that $\|u\|_0 = 1$ and
\begin{equation}\label{eq:weak form of EVP mixed}
\begin{cases}
a(\mathbf p, \mathbf q)+b(\mathbf q, u) =0 &\forall\, \mathbf q \in \mathbf Q, \\
b(\mathbf p, v)=-\lambda (u, v) & \forall\, v \in V,
\end{cases}
\end{equation}
where $a(\mathbf p, \mathbf q) := \int_\Omega \mathbf p \cdot \mathbf q \,\mathrm{d}x$, $b(\mathbf p, v) := -\int_\Omega \operatorname{div}\mathbf p v \,\mathrm{d}x$. The eigenvalue problem \eqref{eq:weak form of EVP mixed} has a sequence of eigenvalues (\cite{MR2652780}):
\begin{equation*}
	0 < \lambda_1 \le \lambda_2 \le \lambda_3 \le \cdots \nearrow +\infty,
\end{equation*}
and the corresponding eigenfunctions:
\begin{equation*}
	(\mathbf p_1, u_1),\, (\mathbf p_2, u_2),\, (\mathbf p_3, u_3), \dots,
\end{equation*}
which  satisfy $(u_i,u_j) = \delta_{ij}$ ($\delta_{ij}$ is the Kronecker delta).

\subsection{HDG finite element discretization}
\label{subs:HDG finite element discretization scheme}

The HDG method yields a scalar approximation $u_h$ to the solution $u$, a vector approximation $\mathbf p_h$ to the flux $\mathbf p$, and a scalar approximation $\hat{u}_h$ to the trace of $u$ on $\mathcal{E}_h$ in the following spaces, respectively:
\begin{align*}
V_h&=\left\{v_h \in L^2(\Omega):\left.v_h\right|_K \in V(K), \quad\forall\, K \in \mathcal{T}_h\right\}, \\
\mathbf Q_h&=\left\{\mathbf q_h \in \mathbf L^2(\Omega):\left.\mathbf q_h\right|_K \in \mathbf Q(K), \quad\forall\, K \in \mathcal{T}_h\right\}, \\
\hat{V}_h&=\left\{\hat{v}_h \in L^2\left(\mathcal{E}_h\right):\left.\hat{v}_h\right|_e \in \hat{V}(e), \quad\forall\, e \in \mathcal{E}_h^i,\quad\left.\hat{v}_h\right|_{\mathcal{E}_h^\partial}=0\right\},
\end{align*}
where $V(K)$, $\mathbf Q(K)$ and $\hat{V}(e)$ are some local spaces. Note that functions in these spaces are not necessarily continuous across element interfaces. For convenience, we define $\widetilde{V}_h:=V_h \times \hat{V}_h$, i.e., for any $\widetilde{v}_{h}\in \widetilde{V}_{h},$ we have $\widetilde{v}_{h}=\{v_{h},\hat{v}_{h}\}$, where $v_{h}\in V_{h}$ and $\hat{v}_{h}\in \hat{V}_{h}$. The discretization variational form of the HDG method for \eqref{eq:EVP} reads: Find $(\lambda_h, \mathbf p_h, \widetilde{u}_h) \in \mathbb{R} \times \mathbf Q_h \times \widetilde{V}_h$, such that $\|u_h\|_0 = 1$ and
\begin{equation}\label{eq:EVP HDG discretization}
\begin{cases}
a_h(\mathbf p_h, \mathbf q_h)+b_h(\mathbf q_h, \widetilde{u}_h)=0, \\
b_h(\mathbf p_h, \widetilde{v}_h)+c_h(\widetilde{u}_h, \widetilde{v}_h)=-\lambda_h(u_h, v_h)_{\mathcal{T}_h} \quad \forall\, (\mathbf q_h, \widetilde{v}_h) \in \mathbf Q_h \times \widetilde{V}_h,
\end{cases}
\end{equation}
where
\begin{equation*}
\begin{cases}
a_h(\mathbf p_h, \mathbf q_h)&:=(\mathbf p_h, \mathbf q_h)_{\mathcal{T}_h}, \\
b_h(\mathbf q_h, \widetilde{u}_h)&:=-(u_h, \operatorname{div} \mathbf q_h)_{\mathcal{T}_h}+\langle\hat{u}_h, \mathbf q_h \cdot \mathbf n_K\rangle_{\partial \mathcal{T}_h} \\
&=(\nabla_h u_h, \mathbf q_h)_{\mathcal{T}_h} - \langle u_h-\hat{u}_h, \mathbf q_h \cdot \mathbf n_K\rangle_{\partial \mathcal{T}_h}, \\
c_h(\widetilde{u}_h, \widetilde{v}_h)&:=-\tau\langle u_h-\hat{u}_h, v_h-\hat{v}_h\rangle_{\partial \mathcal{T}_h},
\end{cases}
\end{equation*}
with $\tau$ being the stabilization parameter.

The HDG method \eqref{eq:EVP HDG discretization} can be written in a compact form: Find $(\lambda_h, \mathbf p_h, \widetilde{u}_h) \in \mathbb{R} \times \mathbf Q_h \times \widetilde{V}_h$, such that $\|u_h\|_0 = 1$ and
\begin{equation}
A_h((\mathbf p_h, \widetilde{u}_h),(\mathbf q_h, \widetilde{v}_h))=-\lambda_h(u_h, v_h)_{\mathcal{T}_h} \quad \forall\, (\mathbf q_h, \widetilde{v}_h) \in \mathbf Q_h \times \widetilde{V}_h,
\end{equation}
where 
\begin{equation*}
A_h((\mathbf p_h, \widetilde{u}_h),(\mathbf q_h, \widetilde{v}_h))=a_h(\mathbf p_h, \mathbf q_h)+b_h(\mathbf q_h, \widetilde{u}_h)+b_h(\mathbf p_h, \widetilde{v}_h)+c_h(\widetilde{u}_h, \widetilde{v}_h) .
\end{equation*}
From \cite{MR606505}, the HDG discretization \eqref{eq:EVP HDG discretization} admits a sequence of discrete eigenvalues:
\begin{equation*}
0 < \lambda_{1,h} \le \lambda_{2,h} \le \cdots \le \lambda_{N,h},
\end{equation*}
and the corresponding eigenfunctions:
\begin{equation*}
(\mathbf p_{1,h},\widetilde{u}_{1,h}),\, (\mathbf p_{2,h},\widetilde{u}_{2,h}), \dots, (\mathbf p_{N,h},\widetilde{u}_{N,h}),
\end{equation*}
which  satisfy $(u_{i,h},u_{j,h})_{\mathcal{T}_h} = \delta_{ij}$ ($N = \mathrm{dim}\  V_h$). 

\section{Gradient-based HDG discretization}
\label{sec:Gradient-based HDG discretization}
The purpose of this section is to investigate upper and lower bound properties of eigenvalues computed with gradient-based HDG method. 
Throughout this section, we are going to choose the following local spaces:
\begin{equation*}
V(K) = \mathcal{P}_{k}(K), \quad \mathbf Q(K) = \bm{\mathcal{P}}_{k-1}(K), \quad \hat{V}(e) = \mathcal{P}_k(e), \quad k \ge 1.
\end{equation*}
Also, we choose $\tau = \gamma h_K^{-1}$ in \eqref{eq:EVP HDG discretization}, where the parameter $\gamma$ is a positive constant on each edge in $\mathcal{E}_h$. Here, $\bm{\mathcal{P}}_{k-1}(K) := [\mathcal{P}_{k-1}(K)]^2$, and $\mathcal{P}_{k}(K)$ is the polynomial space of degree at most $k$ on $K$. For simplicity, in the following we only consider the local spaces with $k = 1$.

For any $\widetilde{v}_h \in \widetilde{V}_h$ and $\mathbf q_h \in \mathbf Q_h$, we may define the following parameter-dependent norms as:
\begin{equation}\label{eq:parameter-dependent norms gradient-based} 
\begin{aligned}
& \left\|\widetilde{v}_h\right\|_{\widetilde{1}, \gamma, h}^2=\left(\nabla_h v_h, \nabla_h v_h\right)_{\mathcal{T}_h} + \gamma \sum_{K \in \mathcal{T}_h} h_K^{-1}\left\langle v_h-\hat{v}_h, v_h-\hat{v}_h\right\rangle_{\partial K}, \\
& \left\|\mathbf q_h\right\|^2=(\mathbf q_h, \mathbf q_h)_{\mathcal{T}_h}.
\end{aligned}
\end{equation}

Before discussing the upper or lower bound properties of eigenvalues, we need to introduce a priori estimates of the corresponding source problem. Consider the following Poisson equation (rewritten as a first order system):
\begin{equation}\label{eq:BVP}
\begin{cases}
\mathbf r+\nabla w=0 & \text { in } \Omega, \\
\operatorname{div} \mathbf r=f & \text { in } \Omega, \\
w = 0 & \text { on } \partial\Omega.
\end{cases}
\end{equation}
where $f \in L^2(\Omega)$. The HDG discretization for \eqref{eq:BVP} reads: Find $(\mathbf r_h, \widetilde{w}_h) \in \mathbf Q_h \times \widetilde{V}_h$, such that
\begin{equation}\label{eq:BVP HDG discretization}
A_h((\mathbf r_h, \widetilde{w}_h),(\mathbf q_h, \widetilde{v}_h))=-(f, v_h)_{\mathcal{T}_h} \quad \forall\, (\mathbf q_h, \widetilde{v}_h) \in \mathbf Q_h \times \widetilde{V}_h.
\end{equation}
For convenience of notations, we denote by $(\mathbf r, w)$, $(\mathbf r_h, w_h)$ the unique solution of \eqref{eq:BVP} and the  finite element solution of \eqref{eq:BVP HDG discretization}, respectively,  and by $(\mathbf p_h, \widetilde{u}_h)$ the finite element solution of eigenvalue problem \eqref{eq:EVP HDG discretization} for distinction in the following context.

The following lemma shows the consistency property of stabilized HDG methods (see \cite{MR3895795}). Note that Lemma \ref{lem:consistency} holds for both two types of HDG methods considered in this paper.
\begin{lmm}\label{lem:consistency}
Let $(\mathbf r, w)$ be the solution of \eqref{eq:BVP}. Then it holds that
\begin{equation}\label{eq:consistency}
A_h((\mathbf r, w), (\mathbf q_h, \widetilde{v}_h)) = -(f,v_h)_{\mathcal{T}_h} \quad \forall\, (\mathbf q_h, \widetilde{v}_h) \in \mathbf Q_h \times \widetilde{V}_h.
\end{equation}
\end{lmm}

To obtain the convergence of the discrete eigenpairs, we need the following basic results.
\begin{thrm}\label{thm:a priori estimate BVP gradient-based HDG}
Let $(\mathbf r, w) \in \mathbf L^2(\Omega) \times H_0^1(\Omega)$ be the true solution of the source problem \eqref{eq:BVP} with $\mathbf r \in \mathbf H^1(\Omega)$, $w \in H^2(\Omega)$, and $(\mathbf r_h, w_h, \hat{w}_h) \in \mathbf Q_h \times V_h \times \hat{V}_h$ be the solution of \eqref{eq:BVP HDG discretization} with $\tau = \gamma h_K^{-1}$. Then for any $\gamma >0$, we have
\begin{equation}\label{eq:gradient-based BVP priori estimate}
\left\|\mathbf r-\mathbf r_h\right\|+\left\|w-\widetilde{w}_h\right\|_{\widetilde{1}, \gamma, h} \leq C_{\gamma,1} h\left(|\mathbf r|_1+|w|_2\right),
\end{equation}
and
\begin{equation}\label{eq:BVP L2 priori estimate}
\|w-w_h\|_0 \le C_{\gamma,2} h^2 |w|_2,
\end{equation}
where $C_{\gamma,1}$, $C_{\gamma,2}$ are constants independent of the mesh size $h$.
\end{thrm}

\begin{proof}
The proof of \eqref{eq:gradient-based BVP priori estimate}, as presented in reference \cite{MR4191130}, can be easily adapted with minor modifications. Moreover, it should be noted that the constant $C$ in the inequality may be dependent on the parameter $\gamma$.

Next, we prove \eqref{eq:BVP L2 priori estimate}. We rewrite the HDG method \eqref{eq:BVP HDG discretization} of the source problem as:
\begin{equation}\label{eq:BVP HDG discretization rewritten gradient-based}
\begin{cases}
(\mathbf r_h, \mathbf q_h)_{\mathcal{T}_h}+(\nabla_h w_h, \mathbf q_h)_{\mathcal{T}_h} - \langle w_h - \hat{w}_h, \mathbf q_h \cdot \mathbf n\rangle_{\partial \mathcal{T}_h}=0  \qquad\qquad \forall\, (\mathbf q_h, v_h, \hat{v}_h) \in \mathbf Q_h \times V_h \times \hat{V}_h,\\
-(\mathbf r_h, \nabla_h v_h)_{\mathcal{T}_h} + \langle\mathbf r_h \cdot \mathbf n, v_h - \hat{v}_h\rangle_{\partial \mathcal{T}_h}+\gamma\langle h_K(w_h-\hat{w}_h), v_h-\hat{v}_h\rangle_{\partial \mathcal{T}_h}=(f, v_h)_{\mathcal{T}_h}.
\end{cases}
\end{equation}
The weak form of the source problem appropriate for the mixed method reads: Find $(\mathbf r, w) \in \mathbf L^2(\Omega) \times H_0^1(\Omega)$, such that
\begin{equation}\label{eq:weak form of the source problem gradient-based}
\begin{cases}
(\mathbf r, \mathbf q) + (\nabla w, \mathbf q) = 0 & \forall\, \mathbf q \in \mathbf L^2(\Omega), \\
-(\mathbf r, \nabla v) = (f,v) & \forall\, v \in H_0^1(\Omega).
\end{cases}
\end{equation}
For any given $\psi \in L^2(\Omega)$, let $\varphi \in H^2(\Omega)\cap H_0^1(\Omega)$ satisfy the following auxiliary problem
\begin{equation}\label{eq:auxiliary problem gradient-based}
-\Delta \varphi=\psi \quad \text { in } \Omega,\left.\quad \varphi\right|_{\partial \Omega}=0.
\end{equation}
By integrating by parts, we see that
\begin{equation*}
(w-w_h, \psi)_{\mathcal{T}_h} = -(w-w_h, \operatorname{div} (\nabla \varphi))_{\mathcal{T}_h} = (\nabla_h (w-w_h),\nabla \varphi)_{\mathcal{T}_h} - \langle w-w_h, \nabla \varphi \cdot \mathbf n\rangle_{\partial \mathcal{T}_h}.
\end{equation*}

Now we choose $\mathbf q_h  = Q_h (\nabla \varphi)$, $v_h = \pi_h \varphi \in V_h \cap H_0^1(\Omega)$, $\mathbf q = \mathbf q_h$, $v = v_h$ in \eqref{eq:BVP HDG discretization rewritten gradient-based} and \eqref{eq:weak form of the source problem gradient-based}, where $Q_h$ denotes $L^2$-projector from $\mathbf L^2(\Omega)$ onto $\mathbf Q_h$, and $\pi_h$ denotes the interpolation of $\varphi$ into the conforming piecewise linear finite element space. Since $v_h \in H_0^1(\Omega)$, we can choose $\hat{v}_h = v_h|_{\mathcal{E}_h}$. Thus subtracting \eqref{eq:BVP HDG discretization rewritten gradient-based} from \eqref{eq:weak form of the source problem gradient-based} leads to
\begin{equation*}
\begin{cases}
(\mathbf r-\mathbf r_h, Q_h (\nabla \varphi))_{\mathcal{T}_h}+(\nabla_h(w - w_h), Q_h (\nabla \varphi))_{\mathcal{T}_h} + \langle w_h - \hat{w}_h, Q_h (\nabla \varphi) \cdot \mathbf n\rangle_{\partial \mathcal{T}_h}=0,\\
-(\mathbf r-\mathbf r_h, \nabla (\pi_h \varphi))_{\mathcal{T}_h}=0.
\end{cases}
\end{equation*}
Then we have
\begin{align}
(w - w_h, \psi)_{\mathcal{T}_h} &= (\nabla_h (w-w_h),\nabla \varphi - Q_h (\nabla \varphi))_{\mathcal{T}_h} - (\mathbf r-\mathbf r_h, Q_h (\nabla \varphi) - \nabla (\pi_h \varphi))_{\mathcal{T}_h} \notag\\
&\quad - \langle w_h - \hat{w}_h, Q_h (\nabla \varphi) \cdot \mathbf n\rangle_{\partial \mathcal{T}_h} - \langle w-w_h, \nabla \varphi \cdot \mathbf n\rangle_{\partial \mathcal{T}_h} \notag\\
&= (\nabla_h (w-w_h),\nabla \varphi - Q_h (\nabla \varphi))_{\mathcal{T}_h} - (\mathbf r-\mathbf r_h, Q_h (\nabla \varphi) - \nabla (\pi_h \varphi))_{\mathcal{T}_h} \notag\\
&\quad - \langle w_h - \hat{w}_h, Q_h (\nabla \varphi) \cdot \mathbf n - \nabla \varphi \cdot \mathbf n\rangle_{\partial \mathcal{T}_h}, \label{eq:(w_wh, psi)}
\end{align}
where $\langle w, \nabla \varphi \cdot \mathbf n\rangle_{\partial \mathcal{T}_h} = \langle \hat{w}_h, \nabla \varphi \cdot \mathbf n\rangle_{\partial \mathcal{T}_h} = 0$ was used in the last step. 
Using the trace theorem, the approximations of $L^2$-projection and interpolation operator $\pi_h$, we obtain
\begin{align}
\|\nabla \varphi - Q_h (\nabla \varphi)\|_0 &\le Ch |\varphi|_2, \notag\\
\|Q_h (\nabla \varphi) - \nabla (\pi_h \varphi)\|_0 &\le \|Q_h (\nabla \varphi) - \nabla \varphi\|_0 + |\varphi -  \pi_h \varphi|_1 \le Ch|\varphi|_{2}, \notag\\
\sum_{K\in \mathcal{T}_h}\|Q_h (\nabla \varphi) \cdot \mathbf n - \nabla \varphi \cdot \mathbf n\|_{0,\partial K}^2 
&\le C\sum_{K\in \mathcal{T}_h}\left(h_K|Q_h (\nabla \varphi) - \nabla \varphi|_{1,K}^2 + h_K^{-1}\|Q_h (\nabla \varphi) - \nabla \varphi\|_{0,K}^2\right) \notag\\
&\le Ch|\varphi|_{2}^{2}+ Ch^{-1}\|Q_h (\nabla \varphi) - \nabla \varphi\|_{0}^2
\le Ch|\varphi|_{2}^{2}. \label{eq:Qh approximation properties}
\end{align}
By \eqref{eq:parameter-dependent norms gradient-based} and \eqref{eq:gradient-based BVP priori estimate}, we get
\begin{equation}\label{eq:wh_hat wh}
\sum_{K\in \mathcal{T}_h}\|w_h - \hat{w}_h\|_{0,\partial K}^2 \le C_{\gamma,1}^2\frac{h^3}{\gamma}.
\end{equation}
Combining \eqref{eq:(w_wh, psi)}, \eqref{eq:Qh approximation properties}, \eqref{eq:wh_hat wh} and \eqref{eq:gradient-based BVP priori estimate} for $\|\nabla_h (w-w_h)\|_0$ and $\|\mathbf r - \mathbf r_h\|$, we have
\begin{align*}
(w - w_h, \psi)_{\mathcal{T}_h} &\le Ch^2 |w|_2 |\varphi|_2 + \left(\sum_{K\in \mathcal{T}_h}\|w_h - \hat{w}_h\|_{0,\partial K}^2\right)^{\frac{1}{2}}\left(\sum_{K\in \mathcal{T}_h}\|Q_h (\nabla \varphi) \cdot \mathbf n - \nabla \varphi \cdot \mathbf n\|_{0,\partial K}^2\right)^{\frac{1}{2}} \\ 
&\le C_{\gamma,2} h^2 |w|_2 |\varphi|_2,
\end{align*}
so that \eqref{eq:BVP L2 priori estimate} follows from elliptic regularity for the Dirichlet problem.
\end{proof}

By virtue of Theorem \ref{thm:a priori estimate BVP gradient-based HDG} and the well-known convergence of the mixed eigenpairs (see \cite{MR606505}), we obtain the convergence of the discrete eigenpairs. The arguments of proof for these consequences are analogous to \cite{MR2220929}, and we shall not repeat them.

\begin{thrm}[Convergence of eigenvalues and eigenfunctions]\label{thm:Convergence of eigenvalues and eigenfunctions gradient-based}
Let $(\lambda_h, \mathbf p_h, \widetilde{u}_h)$ with $\|u_h\|_0 = 1$ be an eigenpair of discrete problem \eqref{eq:EVP HDG discretization}, and choose the local spaces $V(K) = \mathcal{P}_1(K)$, $\mathbf Q(K) = \bm{\mathcal{P}}_0(K)$, $\hat{V}(e) = \mathcal{P}_1(e)$. Then there exists an eigenpair $(\lambda, \mathbf p, u)$ of problem $\eqref{eq:EVP}$ with $\|u\|_0 = 1$ such that
\begin{align}
|\lambda-\lambda_{h}| & \leq C_{\gamma,\lambda} h^2, \\
\|\mathbf p-\mathbf p_h\|_0 & \leq C_{\gamma,\lambda} h, \label{eq:p_ph L2 estimate EVP gradient-based HDG}\\
\|u-u_h\|_0 & \leq C_{\gamma,\lambda} h^2. \label{eq:u_uh L2 estimate EVP gradient-based HDG}
\end{align}
\end{thrm}

The saturation condition is a crucial ingredient for investigating lower or upper bounds of eigenvalues, for which a rigorous proof can be found in \cite{MR3254372,MR3120579}.
\begin{lmm}[Saturation condition]\label{lem:The Saturation Condition}
	Assume $\mathbf p \in \mathbf H^1(\Omega)$, $u \in H^2(\Omega)$, and choose the local spaces $V(K) = \mathcal{P}_1(K)$, $\mathbf Q(K) = \bm{\mathcal{P}}_0(K)$. Then the following lower bounds of the approximation error hold when the partition $\mathcal{T}_h$ is quasi-uniform:
	\begin{align}
	 \inf_{\mathbf q_h \in \mathbf Q_h}\|\mathbf p-\mathbf q_h\|_{0} &\ge Ch, \\
	\inf_{v_h \in V_h \cap H^1(\Omega)}|u-v_h|_{1} &\ge Ch.
	\end{align}
\end{lmm}

\subsection{Upper Bounds for Eigenvalues}
\par In this subsection, we shall derive an identity for discrete eigenpairs by the HDG method, and then achieve upper bounds for eigenvalues with proper penalty parameters. For convenience, let $\mathcal{W} := (\mathbf r,w)$ be the true solution of the source problem \eqref{eq:BVP}, and $\mathcal{U} := (\mathbf p,u)$ be the solution of the eigenproblem \eqref{eq:EVP}. Then \eqref{eq:consistency} can be abbreviated as 
\begin{equation*}
A_h(\mathcal{W},\mathcal{V}_h) = -(f,v_h)_{\mathcal{T}_h} \quad \forall\,\mathcal{V}_h\in \mathbf U_h := \mathbf Q_h \times \widetilde{V}_h,
\end{equation*}
where $\mathcal{V}_h := (\mathbf q_h,\widetilde{v}_h)$.
Take $f = \lambda u$, then $w = u$, $\mathbf r = \mathbf p$ and
\begin{equation}\label{eq:consistency1}
A_h(\mathcal{U}, \mathcal{V}_h) = -\lambda(u,v_h)_{\mathcal{T}_h} \quad \forall\,\mathcal{V}_h \in \mathbf U_h.
\end{equation}
Thus, from \eqref{eq:consistency1}, $A_h(\mathcal{U}, \mathcal{U}) = a_h(-\nabla u, -\nabla u) + 2b_h(-\nabla u,u) = -\lambda$ and $A_h(\mathcal{U}_h, \mathcal{U}_h) = -\lambda_h$ (with $\mathcal{U}_h := (\mathbf p_h,\widetilde{u}_h)$), we have
\begin{align}
&\quad A_h(\mathcal{U}-\mathcal{U}_h, \mathcal{U}-\mathcal{U}_h) + \lambda \|u-u_h\|_0^2 \notag\\
&=A_h(\mathcal{U}, \mathcal{U}) - 2A_h(\mathcal{U},\mathcal{U}_h) + A_h(\mathcal{U}_h,\mathcal{U}_h) + \lambda(u,u) -2\lambda(u,u_h)_{\mathcal{T}_h} + \lambda(u_h,u_h)_{\mathcal{T}_h} \label{eq:lam_lamh equality} \\
&=-\lambda-\lambda_h+2\lambda = \lambda-\lambda_h. \notag
\end{align}

\begin{thrm}\label{thm:Upper Bounds for Eigenvalues gradient-based}
	 Suppose that $\gamma$ is large enough, then $\lambda \le \lambda_h$ for sufficiently small $h$.
\end{thrm}
\begin{proof}
It follows from \eqref{eq:lam_lamh equality} and the definition of $A_h(\cdot,\cdot)$ that
\begin{align*}
\lambda-\lambda_h
&= A_h(\mathcal{U}-\mathcal{U}_h, \mathcal{U}-\mathcal{U}_h) + \lambda \|u-u_h\|_0^2 \\
&=a_h(\mathbf p - \mathbf p_h,\mathbf p - \mathbf p_h) + 2b_h(\mathbf p - \mathbf p_h,u - \widetilde{u}_h) + c_h(u - \widetilde{u}_h,u - \widetilde{u}_h) + \lambda \|u-u_h\|_0^2 \\
&= \|\mathbf p-\mathbf p_h\|_{0}^2 + 2(\nabla_h(u-u_h),\mathbf p-\mathbf p_h)_{\mathcal{T}_h} + 2\langle u_h-\hat{u}_h,(\mathbf p-\mathbf p_h)\cdot \mathbf n\rangle_{\partial \mathcal{T}_h} \\
&\quad - \gamma \sum_{K \in \mathcal{T}_h}h^{-1}\|u_h-\hat{u}_h\|_{0,\partial K}^2 + \lambda \|u-u_h\|_0^2 \\
&= -\|\mathbf p-\mathbf p_h\|_{0}^2 - 2(\mathbf p_h+\nabla_h u_h,\mathbf p-\mathbf p_h)_{\mathcal{T}_h} + 2\langle u_h-\hat{u}_h,(\mathbf p-\mathbf p_h)\cdot \mathbf n\rangle_{\partial \mathcal{T}_h} \\
&\quad - \gamma \sum_{K \in \mathcal{T}_h}h^{-1}\|u_h-\hat{u}_h\|_{0,\partial K}^2 + \lambda \|u-u_h\|_0^2.
\end{align*}

Now we estimate the two terms $\|\mathbf p_h+\nabla_h u_h\|_0$ and $|2\langle u_h-\hat{u}_h,(\mathbf p-\mathbf p_h)\cdot \mathbf n\rangle_{\partial \mathcal{T}_h}|$. First we rewrite the first equation of \eqref{eq:EVP HDG discretization} as:
\begin{equation}\label{eq:ph nabla uh to uh_hat uh}
	(\mathbf p_h+\nabla_h u_h, \mathbf q_h)_{\mathcal{T}_h} = \langle u_h-\hat{u}_h, \mathbf q_h \cdot \mathbf n_K\rangle_{\partial \mathcal{T}_h} \quad \forall\, \mathbf q_h \in \mathbf Q_h.
\end{equation}
From the trace theorem on element $K$, we obtain
\begin{equation*}
	(\mathbf p_h+\nabla_h u_h, \mathbf q_h)_{\mathcal{T}_h} \le \sum_{K \in \mathcal{T}_h}\|u_h-\hat{u}_h\|_{0,\partial K}\|\mathbf q_h\|_{0,\partial K} \le C\sum_{K \in \mathcal{T}_h}h^{-\frac{1}{2}}\|u_h-\hat{u}_h\|_{0,\partial K}\|\mathbf q_h\|_{0,K},
\end{equation*}
which yields
\begin{equation}\label{eq:ph+nabla uh to penalty ineq}
	\|\mathbf p_h+\nabla_h u_h\|_0 \le \left(C_1h^{-1}\sum_{K \in \mathcal{T}_h}\|u_h-\hat{u}_h\|_{0,\partial K}^2\right)^{\frac{1}{2}}.
\end{equation}
Then, by the Cauchy-Schwarz inequality, the trace theorem and Lemma \ref{lem:The Saturation Condition}, we have
\begin{align}
&\quad 2\langle u_h-\hat{u}_h, (\mathbf p-\mathbf p_h) \cdot \mathbf n_K\rangle_{\partial \mathcal{T}_h} \notag\\
&\le \sum_{K \in \mathcal{T}_h}\left(\varepsilon h \|\mathbf p-\mathbf p_h\|_{0,\partial K}^2 + \varepsilon^{-1}h^{-1}\|u_h-\hat{u}_h\|_{0,\partial K}^2\right) \notag\\
&\le C\varepsilon\sum_{K \in \mathcal{T}_h}\left(h^2 |\mathbf p|_{1,K}^2 + \|\mathbf p-\mathbf p_h\|_{0,K}^2\right) + \varepsilon^{-1}\sum_{K \in \mathcal{T}_h}h^{-1}\|u_h-\hat{u}_h\|_{0,\partial K}^2 \notag\\
&\le C_2\varepsilon\|\mathbf p-\mathbf p_h\|_{0}^2 + \varepsilon^{-1}\sum_{K \in \mathcal{T}_h}h^{-1}\|u_h-\hat{u}_h\|_{0,\partial K}^2, \label{eq:Young ineq}
\end{align}
where $\varepsilon > 0$ is arbitrary. Finally, by choosing
\begin{equation*}
	\varepsilon=\frac{1}{4 C_2} \quad \text { and } \quad \gamma \ge \frac{1}{2}+4(C_1+C_2),
\end{equation*}
we have, in view of \eqref{eq:ph+nabla uh to penalty ineq} and \eqref{eq:Young ineq},
\begin{align*}
\lambda-\lambda_h &\le -\|\mathbf p-\mathbf p_h\|_{0}^2 + 4C_1\sum_{K \in \mathcal{T}_h}h^{-1}\|u_h-\hat{u}_h\|_{0,\partial K}^2 + \frac{1}{4}\|\mathbf p-\mathbf p_h\|_{0}^2 \\
&\quad+ C_2\varepsilon\|\mathbf p-\mathbf p_h\|_{0}^2 + \varepsilon^{-1}\sum_{K \in \mathcal{T}_h}h^{-1}\|u_h-\hat{u}_h\|_{0,\partial K}^2 - \gamma \sum_{K \in \mathcal{T}_h}h^{-1}\|u_h-\hat{u}_h\|_{0,\partial K}^2 + \lambda \|u-u_h\|_0^2 \\
&\le -\frac{1}{2}\|\mathbf p-\mathbf p_h\|_{0}^2 - \frac{1}{2} \sum_{K \in \mathcal{T}_h}h^{-1}\|u_h-\hat{u}_h\|_{0,\partial K}^2 + \lambda \|u-u_h\|_0^2.
\end{align*}
Recall from the convergence of eigenfunctions \eqref{eq:p_ph L2 estimate EVP gradient-based HDG} and \eqref{eq:u_uh L2 estimate EVP gradient-based HDG} that the term $\lambda \|u-u_h\|_0^2$ is a higher order term. This implies $\lambda - \lambda_h$ is negative, regardless of how the a priori estimate of $\|u-u_h\|_0$ depends on stabilization parameters. This concludes the proof.
\end{proof}

\subsection{Lower Bounds for Eigenvalues}
In this subsection, we shall discuss the lower bound for exact eigenvalues.

Define an element-piecewise operator $Q_h^1:L^{2}(\Omega)\to V_h$, i.e., for each $K \in \mathcal{T}_h$, $Q_h^1|_K:L^{2}(K)\to P_1(K)$ satisfies that for any $u\in L^{2}(K)$,
\begin{equation*}
(Q_h^1|_K u,v)_K=(u,v)_K\quad \forall\,v\in P_1(K).
\end{equation*}
Also, define $P_h: H_0^1(\Omega) \rightarrow V_h \cap H_0^1(\Omega)$ to be a Ritz projector, i.e., for any $u\in H_{0}^{1}(\Omega)$,
\begin{equation*}
(\nabla P_{h}u,\nabla v_{h})_{\mathcal{T}_h}=(\nabla u,\nabla v_{h})_{\mathcal{T}_h}\quad \forall\, v_{h}\in V_h \cap H_0^1(\Omega).
\end{equation*}
Suppose $(\lambda, \mathbf p, u)$ to be an eigenpair of problem $\eqref{eq:EVP}$. Let $\mathcal{X}_h:=(\mathbf s_h, \widetilde{x}_h) \in \mathbf U_h$ with $\widetilde{x}_{h}:=\{x_{h},\hat{x}_{h}\}$, where $x_h = Q_h^1 u \in V_h$, $\hat{x}_h = \left.(P_h u)\right|_{\mathcal{E}_h} \in \hat{V}_h$ and $\mathbf s_h = -\nabla P_h u$. It holds that
\begin{align}
\lambda - \lambda_h &= -A_h(\mathcal{U},\mathcal{U}) + A_h(\mathcal{U}_h,\mathcal{U}_h) + 2A_h(\mathcal{U}_h,\mathcal{X}_h) + 2\lambda_h(u_h,x_h)_{\mathcal{T}_h} \notag \\
&= -A_h(\mathcal{U},\mathcal{U}) + A_h(\mathcal{X}_h,\mathcal{X}_h) - A_h(\mathcal{U}_h,\mathcal{U}_h) + 2A_h(\mathcal{U}_h,\mathcal{X}_h) - A_h(\mathcal{X}_h,\mathcal{X}_h) \notag \\
&\quad + 2A_h(\mathcal{U}_h,\mathcal{U}_h) + 2\lambda_h(u_h,u)_{\mathcal{T}_h} \notag \\
&= -A_h(\mathcal{U},\mathcal{U}) + A_h(\mathcal{X}_h,\mathcal{X}_h) - A_h(\mathcal{U}_h-\mathcal{X}_h,\mathcal{U}_h-\mathcal{X}_h) - \lambda_h\|u-u_h\|_0^2. \label{eq:lam_lamh equality for lower bounds}
\end{align}
\begin{lmm}[{see \cite[Lemma 3.16]{MR4476865}}]\label{lem:xh_hat xh priori estimate}
Suppose $(\lambda, \mathbf p, u)$ to be an eigenpair of problem $\eqref{eq:EVP}$, and let $x_h = Q_h^1 u$, $\hat{x}_h = \left.(P_h u)\right|_{\mathcal{E}_h}$. Then it holds that
\begin{equation*}
\sum_{K \in \mathcal{T}_h}h_K^{-1}\|x_h-\hat{x}_h\|_{0,\partial K}^2 \le C_1 h^2.
\end{equation*}
\end{lmm}

Now we are ready to state the lower bound properties of eigenvalues.
\begin{thrm}\label{thm:Lower Bounds for Eigenvalues gradient-based}
	Suppose that $\gamma$ is small enough, then $\lambda \ge \lambda_h$ for sufficiently small $h$.
\end{thrm}

\begin{proof}
We examine the two terms $-A_h(\mathcal{U},\mathcal{U}) + A_h(\mathcal{X}_h,\mathcal{X}_h)$ and $- A_h(\mathcal{U}_h-\mathcal{X}_h,\mathcal{U}_h-\mathcal{X}_h)$ independently. First, by integrating by parts and the definition of $\hat{x}_h$ and $\mathbf s_h$, we have 
\begin{align}
\langle x_h-\hat{x}_h, \mathbf s_h\cdot \mathbf n\rangle_{\partial \mathcal{T}_h} &= \langle x_h-P_h u, \mathbf s_h\cdot \mathbf n\rangle_{\partial \mathcal{T}_h} \notag\\
&= (\nabla_h x_h-\nabla P_h u,\mathbf s_h)_{\mathcal{T}_h} \notag\\
&= (\nabla_h x_h+\mathbf s_h,\mathbf s_h)_{\mathcal{T}_h},\label{eq:xh_hat xh Green formula}
\end{align}
from which it follows that
\begin{align}
&\quad -A_h(\mathcal{U},\mathcal{U}) + A_h(\mathcal{X}_h,\mathcal{X}_h) \notag \\
&= |u|_1^2 + (\mathbf s_h,\mathbf s_h)_{\mathcal{T}_h} + 2(\nabla_h x_h, \mathbf s_h)_{\mathcal{T}_h} - 2\langle x_h-\hat{x}_h, \mathbf s_h\cdot \mathbf n\rangle_{\partial \mathcal{T}_h} - \gamma h^{-1}\|x_h-\hat{x}_h\|_{\partial \mathcal{T}_h}^2 \notag \\
&= |u|_1^2 + (\mathbf s_h,\mathbf s_h)_{\mathcal{T}_h} + 2(\nabla_h x_h, \mathbf s_h)_{\mathcal{T}_h} - 2(\nabla_h x_h+\mathbf s_h,\mathbf s_h)_{\mathcal{T}_h}- \gamma h^{-1}\|x_h-\hat{x}_h\|_{\partial \mathcal{T}_h}^2 \notag \\
&= |u|_1^2 - \|\mathbf s_h\|_0^2- \gamma h^{-1}\|x_h-\hat{x}_h\|_{\partial \mathcal{T}_h}^2 \notag \\
&= |u-P_h u|_1^2- \gamma h^{-1}\|x_h-\hat{x}_h\|_{\partial \mathcal{T}_h}^2, \label{eq:the 1st term in equality for lower bounds}
\end{align}
where the property of orthogonal projector $P_h$ was used in the last step. Then from \eqref{eq:ph nabla uh to uh_hat uh}, the definition of $A_h(\cdot,\cdot)$ and the same technique as \eqref{eq:xh_hat xh Green formula}, we have
\begin{align}
&\quad- A_h(\mathcal{U}_h-\mathcal{X}_h,\mathcal{U}_h-\mathcal{X}_h) \notag \\
&= -\|\mathbf p_h-\mathbf s_h\|_0^2-2(\nabla_h(u_h-x_h),\mathbf p_h-\mathbf s_h)_{\mathcal{T}_h}+2\langle u_h-\hat{u}_h-(x_h-\hat{x}_h), (\mathbf p_h-\mathbf s_h)\cdot \mathbf n\rangle_{\partial \mathcal{T}_h} \notag \\
&\quad +\gamma h^{-1}\|u_h-\hat{u}_h-(x_h-\hat{x}_h)\|_{\partial \mathcal{T}_h}^2 \notag \\
&= -\|\mathbf p_h-\mathbf s_h\|_0^2-2(\nabla_h(u_h-x_h),\mathbf p_h-\mathbf s_h)_{\mathcal{T}_h} +2(\mathbf p_h+\nabla_h u_h,\mathbf p_h-\mathbf s_h)_{\mathcal{T}_h} \notag \\
&\quad- 2(\nabla_h x_h + \mathbf s_h,\mathbf p_h-\mathbf s_h)_{\mathcal{T}_h} + \gamma h^{-1}\|u_h-\hat{u}_h-(x_h-\hat{x}_h)\|_{\partial \mathcal{T}_h}^2 \notag \\
&= \|\mathbf p_h-\mathbf s_h\|_0^2 + \gamma h^{-1}\|u_h-\hat{u}_h-(x_h-\hat{x}_h)\|_{\partial \mathcal{T}_h}^2>0. \label{eq:the 3rd term in equality for lower bounds}
\end{align}
Finally combining \eqref{eq:lam_lamh equality for lower bounds}, \eqref{eq:the 1st term in equality for lower bounds}, \eqref{eq:the 3rd term in equality for lower bounds}, Lemma \ref{lem:The Saturation Condition}, Lemma \ref{lem:xh_hat xh priori estimate} and a priori estimate \eqref{eq:u_uh L2 estimate EVP gradient-based HDG}, the following holds
\begin{align*}
\lambda - \lambda_h &= |u-P_h u|_1^2- \gamma h^{-1}\|x_h-\hat{x}_h\|_{\partial \mathcal{T}_h}^2\\
&\quad + \|\mathbf p_h-\mathbf s_h\|_0^2 + \gamma h^{-1}\|u_h-\hat{u}_h-(x_h-\hat{x}_h)\|_{\partial \mathcal{T}_h}^2 - \lambda_h\|u-u_h\|_0^2 \\
&\ge |u-P_h u|_1^2- \gamma h^{-1}\|x_h-\hat{x}_h\|_{\partial \mathcal{T}_h}^2 - \lambda_h\|u-u_h\|_0^2\\
&\ge (C-C_1\gamma)h^2-C_{\gamma, \lambda}h^4,
\end{align*}
which completes the proof.
\end{proof}

\section{Divergence-based HDG discretization}
\label{sec:Divergence-based HDG discretization}
This section is devoted to the divergence-based HDG method with the local spaces:
\begin{equation*}
V(K) = \mathcal{P}_0(K), \quad \mathbf Q(K) = \bm{\mathcal{P}}_1(K), \quad \hat{V}(e) = \mathcal{P}_1(e),
\end{equation*}
and $\tau = \gamma h_K$ in \eqref{eq:EVP HDG discretization}.

Define an interpolation operator $\Pi_h:\mathbf H(\operatorname{div}, \Omega)\to \mathbf Q_h \cap \mathbf H(\operatorname{div}, \Omega)$ and $Q_h^0:L^2(\Omega)\to V_h$ satisfying the following commuting diagram property (see \cite{MR1115205}):
\begin{figure} 
\centering
\begin{tikzpicture}
  \node (A) at (0,0) {$\mathbf Q_h \cap \mathbf H(\operatorname{div}, \Omega)$};
  \node (B) at (3.5,0) {$V_h$};
  \node (C) at (0,1.8) {$\mathbf H(\operatorname{div}, \Omega)$};
  \node (D) at (3.5,1.8) {$L^2(\Omega)$};
  \node (E) at (1.3,1.8) {};
  \node (F) at (3,1.8) {};
  \node (G) at (1.3,0) {};
  \node (H) at (3,0) {};
  
  \draw[->] (G) -- (H) node[midway, above] {$\operatorname{div}$};
  \draw[->] (E) -- (F) node[midway, above] {$\operatorname{div}$};
  \draw[->] (C) -- (A) node[midway, left] {$\Pi_h$};
  \draw[->] (D) -- (B) node[midway, right] {$Q_h^0$};
\end{tikzpicture}
\end{figure}
\noindent where $Q_h^0$ denotes the $L^2$-projection onto $V_h$. The following approximation properties of the interpolation operator $\Pi_h$ and the $L^2$-projection $Q_h^0$ are proved in \cite{MR799685} and references therein:
\begin{align}
\|\mathbf q - \Pi_h \mathbf q\|_0 &\le Ch^{2}|\mathbf q|_{2}, \label{eq:Pi_h L2 estimate} \\
\|\operatorname{div}(\mathbf q - \Pi_h \mathbf q)\|_0 &\le Ch |\operatorname{div}\mathbf q|_1, \\
\|v - Q_h^0 v\|_0 &\le Ch |v|_1, \label{eq: Q0 L2 estimate}
\end{align}
for any $\mathbf q \in \mathbf H^{2}(\Omega)$, $\operatorname{div}\mathbf q \in H^1(\Omega)$, and $v \in H^1(\Omega)$.

Let $Q_{h,e}^1$ be the projection operator defined locally by $L^2(e)$-projection onto $\mathcal{P}_1(e)$ for any $e \in \mathcal{E}_h$. Then we have the following approximation property:
\begin{equation}\label{eq:Qe^1 L2 estimate}
\left\|v - Q_{h,e}^1 v\right\|_{0,e} \le C h_K^{3/2} |v|_{2,K},
\end{equation}
for any $v \in H^2(K)$ with edge $e \subset \partial K$. In fact, it follows from the well-known interpolation error estimate \cite{MR2373954} and the trace theorem that (recall the definition of $\pi_h$ in Section \ref{sec:Gradient-based HDG discretization})
\begin{align*}
\left\|v - Q_{h,e}^1 v\right\|_{0,e} 
\le \left\|v - \pi_h v\right\|_{0,e}
\le C\left(h_K^{1/2} \left|v - \pi_h v\right|_{1,K} + h_K^{-1/2} \left\|v - \pi_h v\right\|_{0,K}\right) \le C h_K^{3/2} |v|_{2,K}.
\end{align*}

For any $\widetilde{v}_h \in \widetilde{V}_h$ and $\mathbf q_h \in \mathbf Q_h$, if $\gamma > 0$, then we may define the following parameter-dependent norms as:
\begin{align*}
& \|\widetilde{v}_h\|_{0, \gamma, h}^2=(v_h, v_h)_{\mathcal{T}_h}+\gamma \sum_{e \in \mathcal{E}_h^i} h_e\langle\hat{v}_h, \hat{v}_h\rangle_e \\
& \|\mathbf q_h\|_{\operatorname{div}, \gamma, h}^2=(\mathbf q_h, \mathbf q_h)_{\mathcal{T}_h}+(\operatorname{div} \mathbf q_h, \operatorname{div} \mathbf q_h)_{\mathcal{T}_h}+\gamma^{-1} \sum_{e \in \mathcal{E}_h^i} h_e^{-1}\left\langle Q_{h,e}^1([\mathbf q_h]), Q_{h,e}^1([\mathbf q_h])\right\rangle_e,
\end{align*}
and if $\gamma = 0$, then
\begin{align*}
& \|\widetilde{v}_h\|_{0, 1, h}^2=(v_h, v_h)_{\mathcal{T}_h}+ \sum_{e \in \mathcal{E}_h^i} h_e\langle\hat{v}_h, \hat{v}_h\rangle_e \\
& \|\mathbf q_h\|_{\operatorname{div}, 1, h}^2=(\mathbf q_h, \mathbf q_h)_{\mathcal{T}_h}+(\operatorname{div} \mathbf q_h, \operatorname{div} \mathbf q_h)_{\mathcal{T}_h}+ \sum_{e \in \mathcal{E}_h^i} h_e^{-1}\left\langle Q_{h,e}^1([\mathbf q_h]), Q_{h,e}^1([\mathbf q_h])\right\rangle_e.
\end{align*}

Similar to Section \ref{sec:Gradient-based HDG discretization}, we first introduce a priori estimates of the corresponding source problem \eqref{eq:BVP}. For any given $\tau = \gamma h_K$, we rewrite the HDG method \eqref{eq:BVP HDG discretization} of the source problem as: Find $(\mathbf r_h, w_h, \hat{w}_h) \in \mathbf Q_h \times V_h \times \hat{V}_h$, such that for any $(\mathbf q_h, v_h, \hat{v}_h) \in \mathbf Q_h \times V_h \times \hat{V}_h$, we have
\begin{equation}\label{eq:BVP HDG discretization rewritten}
\begin{cases}
(\mathbf r_h, \mathbf q_h)_{\mathcal{T}_h}-(w_h, \operatorname{div} \mathbf q_h)_{\mathcal{T}_h}+\langle\hat{w}_h, \mathbf q_h \cdot \mathbf n\rangle_{\partial \mathcal{T}_h}=0,\\
-(\operatorname{div} \mathbf r_h, v_h)_{\mathcal{T}_h}+\langle\mathbf r_h \cdot \mathbf n, \hat{v}_h\rangle_{\partial \mathcal{T}_h}-\gamma\langle h_K(w_h-\hat{w}_h), v_h-\hat{v}_h\rangle_{\partial \mathcal{T}_h}=-(f, v_h)_{\mathcal{T}_h}.
\end{cases}
\end{equation}

The following estimate is proved in \cite{MR799685}.
\begin{lmm}\label{lem:edge estimate lem}
If $(\mathbf r_h, \widetilde{w}_h) \in \mathbf Q_h \times \widetilde{V}_h$ is the solution of discrete problem \eqref{eq:BVP HDG discretization rewritten}, then
\begin{equation}\label{eq:edge estimate}
\|Q_{h,e}^1 w-\hat{w}_h\|_{0,e} \le C\left(h_K^{1/2}\|\mathbf r-\mathbf r_h\|_{0,K} + h_K^{-1/2}\|Q_h^0 w-w_h\|_{0,K}\right).
\end{equation}
\end{lmm}
\begin{rmrk}
Note that Lemma \ref{lem:edge estimate lem} remains valid when considering eigenvalue problems, as the proof of \eqref{eq:edge estimate} relies solely on the form of the first equation in \eqref{eq:BVP HDG discretization rewritten}.
\end{rmrk}

\begin{thrm}\label{thm:a priori estimate BVP divergence-based HDG}
Let $(\mathbf r, w) \in \mathbf H(\operatorname{div}, \Omega) \times L^2(\Omega)$ be the solution of \eqref{eq:BVP} with $\mathbf r \in \mathbf H^2(\Omega)$, $\operatorname{div} \mathbf r \in H^1(\Omega)$, $w \in H^1(\Omega)$, and $(\mathbf r_h, w_h, \hat{w}_h) \in \mathbf Q_h \times V_h \times \hat{V}_h$ be the solution of \eqref{eq:BVP HDG discretization} with $\tau = \gamma h_K$. Then for any $\gamma \ge 0$, we have
\begin{align}
\|\mathbf r-\mathbf r_h\|_{\operatorname{div}, \gamma, h}+\left\|w-\widetilde{w}_h\right\|_{0, \gamma, h} &\leq C_{\gamma,3} h\left(|\mathbf r|_1+|\operatorname{div} \mathbf r|_1+|w|_1\right), \label{eq:divergence-based BVP priori estimate} \\
\|w_h - Q_h^0 w\|_0 &\le C_{\gamma,4}h^2, \label{eq:BVP superconvergence result} \\
\|\mathbf r-\mathbf r_h\|_0 &\le C_{\gamma,5}h^{3/2}. \label{eq:BVP r_r_h L2 estimate}
\end{align}
\end{thrm}
\begin{proof}
The proof of \eqref{eq:divergence-based BVP priori estimate} can be adapted from \cite{MR4191130} with minor modifications and the constant $C$ in the inequality may depend on the parameter $\gamma$. Thus, it suffices to prove \eqref{eq:BVP superconvergence result} and \eqref{eq:BVP r_r_h L2 estimate}.

Note that the discretization \eqref{eq:BVP HDG discretization} is consistent. Subtracting \eqref{eq:BVP HDG discretization rewritten} from \eqref{eq:consistency} leads to the error equations
\begin{equation}\label{eq:BVP error equation}
\begin{cases}
(\mathbf r-\mathbf r_h, \mathbf q_h)_{\mathcal{T}_h}-(Q_h^0 w-w_h, \operatorname{div} \mathbf q_h)_{\mathcal{T}_h}+\langle w-\hat{w}_h, \mathbf q_h \cdot \mathbf n\rangle_{\partial \mathcal{T}_h}=0,\\
-(\operatorname{div} (\mathbf r-\mathbf r_h), v_h)_{\mathcal{T}_h}+\langle (\mathbf r-\mathbf r_h) \cdot \mathbf n, \hat{v}_h\rangle_{\partial \mathcal{T}_h}+\gamma\langle h_K(w_h-\hat{w}_h), v_h-\hat{v}_h\rangle_{\partial \mathcal{T}_h}=0.
\end{cases}
\end{equation}
Let $\psi \in L^2(\Omega)$, and $\varphi \in H^2(\Omega)\cap H_0^1(\Omega)$ satisfy the following auxiliary problem
\begin{equation}\label{eq:auxiliary problem}
-\Delta \varphi=\psi \quad \text { in } \Omega,\left.\quad \varphi\right|_{\partial \Omega}=0.
\end{equation}
By the property of the interpolation operator $\Pi_h$, the projection operator $Q_{h,e}^1$, the first equation of \eqref{eq:BVP error equation} and integration by parts, we obtain
\begin{align}\label{eq:superconvergence eq1}
\begin{aligned}
(Q_h^0 w - w_h, \psi)_{\mathcal{T}_h} 
&= -\left(Q_h^0 w - w_h, Q_{h}^{0}\operatorname{div}\nabla \varphi\right)_{\mathcal{T}_h} \\
&= -\left(Q_h^0 w - w_h, \operatorname{div}\left(\Pi_h \nabla \varphi\right)\right)_{\mathcal{T}_h} \\
&= -(\mathbf r - \mathbf r_h, \Pi_h \nabla \varphi)_{\mathcal{T}_h} \\
&= (\mathbf r - \mathbf r_h, \nabla \varphi - \Pi_h \nabla \varphi)_{\mathcal{T}_h} + (\operatorname{div} (\mathbf r - \mathbf r_h), \varphi)_{\mathcal{T}_h} - \left\langle (\mathbf r - \mathbf r_h) \cdot \mathbf n, Q_{h,e}^1\varphi\right\rangle_{\partial \mathcal{T}_h}.
\end{aligned}
\end{align}
Now we choose $v_h = Q_h^0 \varphi \in V_h$ and $\hat{v}_h = Q_{h,e}^1 \varphi \in \hat{V}_h$ in \eqref{eq:BVP error equation}. By \eqref{eq:superconvergence eq1} and \eqref{eq:BVP error equation}, we get
\begin{align}
(Q_h^0 w - w_h, \psi)_{\mathcal{T}_h} 
&= \left(\mathbf r - \mathbf r_h, \nabla \varphi - \Pi_h \nabla \varphi\right)_{\mathcal{T}_h} + \left(\operatorname{div} (\mathbf r - \mathbf r_h), \varphi - Q_h^0 \varphi\right)_{\mathcal{T}_h} \notag\\
&\quad\,\gamma\left\langle h_K(w_h-\hat{w}_h), Q_h^0 \varphi-Q_{h,e}^1 \varphi\right\rangle_{\partial \mathcal{T}_h} \label{eq:expand item} \\
&:= I_1 + I_2 + I_3. \notag
\end{align}
For the first term on the right hand side of \eqref{eq:expand item}, by \eqref{eq:divergence-based BVP priori estimate} and \eqref{eq:Pi_h L2 estimate}, we have
\begin{equation*}
|I_1| \le C_{\gamma}h^2 \left(|\mathbf r|_1+|\operatorname{div} \mathbf r|_1+|w|_1\right)|\nabla \varphi|_1.
\end{equation*}
By \eqref{eq:divergence-based BVP priori estimate} and \eqref{eq: Q0 L2 estimate}, we have
\begin{equation*}
|I_2| \le C_{\gamma}h^2\left(|\mathbf r|_1+|\operatorname{div} \mathbf r|_1+|w|_1\right)|\varphi|_1.
\end{equation*}
The definition of $\left\|\cdot\right\|_{0, \gamma, h}$ and the a priori estimate \eqref{eq:divergence-based BVP priori estimate} lead to
\begin{equation}\label{eq:w_hat wh}
\sum_{e\in \mathcal{E}_h^i}\|w-\hat{w}_h\|_{0,e}^2 \le C_{\gamma,3}^2\frac{h}{\gamma}.
\end{equation}
By the trace theorem, \eqref{eq: Q0 L2 estimate}, \eqref{eq:Qe^1 L2 estimate}, \eqref{eq:divergence-based BVP priori estimate}, \eqref{eq:w_hat wh} and note that $w|_{\partial \Omega} = \hat{w}_h|_{\partial \Omega} = 0$, we have
\begin{align}
\sum_{K\in \mathcal{T}_h}\|w_h-\hat{w}_h\|_{0,\partial K}^2 &\le 2\sum_{K\in \mathcal{T}_h}\left(\|w_h-w\|_{0,\partial K}^2 + \|w-\hat{w}_h\|_{0,\partial K}^2\right) \notag\\
&\le C\sum_{K\in\mathcal{T}_h}\left(h_K|w|_{1,K}^2 + h_K^{-1}\|w_h-w\|_{0,K}^2\right) + 4\sum_{e\in \mathcal{E}_h^i}\|w-\hat{w}_h\|_{0,e}^2 
\le C_{\gamma} h, \label{eq:wh_hat_wh L2e estimate}
\end{align}
and
\begin{align*}
\sum_{e\in \mathcal{E}_h}\|Q_h^0 \varphi-Q_{h,e}^1 \varphi\|_{0,e}^2 &\le 2\left(\sum_{e\in \mathcal{E}_h}\|Q_h^0 \varphi-\varphi\|_{0,e}^2 + \sum_{e\in \mathcal{E}_h^i}\|\varphi-Q_{h,e}^1 \varphi\|_{0,e}^2\right) \\
&\le C\sum_{K\in \mathcal{T}_h}\left(h_K|\varphi|_{1,K}^2 + h_K^{-1}\|Q_h^0 \varphi-\varphi\|_{0,K}^2\right) + Ch^3|\varphi|_{2}^2 \le Ch\|\varphi\|_{2}^2.
\end{align*}
Combining these inequalities shows that
\begin{equation*}
|I_3| \le C_{\gamma} h^2\|\varphi\|_{2},
\end{equation*}
so that \eqref{eq:BVP superconvergence result} follows from the above inequalities and the assumed elliptic regularity for the Dirichlet problem \eqref{eq:auxiliary problem}.

Next, let us bound $\|\mathbf r - \mathbf r_h\|_0$. In \eqref{eq:BVP error equation}, take $\mathbf q_h = \Pi_h \mathbf r - \mathbf r_h$, $v_h = Q_h^0 w - w_h$ and $\hat{v}_h = Q_{h,e}^1 w - w_h$, we have
\begin{align*}
\|\Pi_h \mathbf r - \mathbf r_h\|_0^2 &= \left(\Pi_h \mathbf r - \mathbf r,\Pi_h \mathbf r - \mathbf r_h\right)_{\mathcal{T}_h} + \left(\mathbf r - \mathbf r_h,\Pi_h \mathbf r - \mathbf r_h\right)_{\mathcal{T}_h} \\
&= \left(\Pi_h \mathbf r - \mathbf r,\Pi_h \mathbf r - \mathbf r_h\right)_{\mathcal{T}_h} + \left(Q_h^0 w-w_h, \operatorname{div} (\Pi_h \mathbf r - \mathbf r_h)\right)_{\mathcal{T}_h} \\
&\quad-\left\langle w-\hat{w}_h, (\Pi_h \mathbf r - \mathbf r_h) \cdot \mathbf n\right\rangle_{\partial \mathcal{T}_h} \\
&= \left(\Pi_h \mathbf r - \mathbf r,\Pi_h \mathbf r - \mathbf r_h\right)_{\mathcal{T}_h} + \left(Q_h^0 w-w_h, \operatorname{div} (\mathbf r - \mathbf r_h)\right)_{\mathcal{T}_h} \\
&\quad-\left\langle Q_{h,e}^1 w-\hat{w}_h, (\mathbf r - \mathbf r_h) \cdot \mathbf n\right\rangle_{\partial \mathcal{T}_h} \\
&= \left(\Pi_h \mathbf r - \mathbf r,\Pi_h \mathbf r - \mathbf r_h\right)_{\mathcal{T}_h} + \gamma\left\langle h_K(w_h-\hat{w}_h), Q_h^0 w - w_h - \left(Q_{h,e}^1 w -\hat{w}_h\right)\right\rangle_{\partial \mathcal{T}_h},
\end{align*}
where we have used the following equalities
\begin{equation*}
\left\langle w-\hat{w}_h, (\Pi_h \mathbf r - \mathbf r_h) \cdot \mathbf n\right\rangle_{\partial \mathcal{T}_h} = \left\langle Q_{h,e}^1 w-\hat{w}_h, (\Pi_h\mathbf r - \mathbf r_h) \cdot \mathbf n\right\rangle_{\partial \mathcal{T}_h} = \left\langle Q_{h,e}^1 w-\hat{w}_h, (\mathbf r - \mathbf r_h) \cdot \mathbf n\right\rangle_{\partial \mathcal{T}_h},
\end{equation*}
because $(\Pi_h \mathbf r - \mathbf r_h) \cdot \mathbf n|_e \in \mathcal{P}_1(e)$ and $(Q_{h,e}^1 w-\hat{w}_h)|_e \in \mathcal{P}_1(e)$ for all $e\in \mathcal{E}_h$, plus the property of projector $Q_{h,e}^1$ and interpolation operator $\Pi_h$.
By the trace theorem, \eqref{eq:BVP superconvergence result}, \eqref{eq:wh_hat_wh L2e estimate} and the Cauchy-Schwarz inequality, we have
\begin{equation*}
\gamma\left\langle h_K(w_h-\hat{w}_h), Q_h^0 w - w_h\right\rangle_{\partial \mathcal{T}_h} \le Ch^{3/2}\sum_{K\in \mathcal{T}_h}h_K^{-1/2}\|Q_h^0 w - w_h\|_{0,K} \le C_{\gamma} h^3.
\end{equation*}
By Lemma \ref{lem:edge estimate lem}, \eqref{eq:divergence-based BVP priori estimate}, \eqref{eq:BVP superconvergence result}, \eqref{eq:wh_hat_wh L2e estimate}, we have
\begin{equation*}
\gamma\left\langle h_K(w_h-\hat{w}_h), Q_{h,e}^1 w -\hat{w}_h\right\rangle_{\partial \mathcal{T}_h} \le Ch^{3/2}\sum_{K\in \mathcal{T}_h} \left(h_K^{1/2}\|\mathbf r-\mathbf r_h\|_{0,K} + h_K^{-1/2}\|Q_h^0 w-w_h\|_{0,K}\right) \le C_{\gamma} h^3.
\end{equation*}
Applying \eqref{eq:Pi_h L2 estimate}, we see that
\begin{equation*}
\|\Pi_h \mathbf r - \mathbf r_h\|_0^2 \le c_1h^2 \|\Pi_h \mathbf r - \mathbf r_h\|_0 + c_2h^3.
\end{equation*}
Then consider $\|\Pi_h \mathbf r - \mathbf r_h\|_0$ as the variable of a quadratic inequality and a fundamental calculation implies that
\begin{equation*}
\|\Pi_h \mathbf r - \mathbf r_h\|_0 \le \frac{1}{2} \left(c_1 h^2 + \sqrt{c_1^2 h^4 + 4c_2 h^3}\right) \le C_{\gamma} h^{3/2},
\end{equation*}
and \eqref{eq:BVP r_r_h L2 estimate} follows from \eqref{eq:Pi_h L2 estimate} and the triangle inequality.
\end{proof}

Let us turn to the analysis of eigenvalues and eigenfunctions. The theorem below on the convergence of the HDG method follows from Theorem \ref{thm:a priori estimate BVP divergence-based HDG} and the classical spectral approximation theory \cite{MR2652780,MR606505}.

\begin{thrm}[Convergence of eigenvalues and eigenfunctions]\label{thm:Convergence of eigenvalues and eigenfunctions divergence-based HDG}
Let $(\lambda_h, \mathbf p_h, \widetilde{u}_h)$ with $\|u_h\|_0 = 1$ be an eigenpair of discrete problem \eqref{eq:EVP HDG discretization}, and choose the local spaces $V(K) = \mathcal{P}_0(K)$, $\mathbf Q(K) = \bm{\mathcal{P}}_1(K)$, $\hat{V}(e) = \mathcal{P}_1(e)$. Then there exists an eigenpair $(\lambda, \mathbf p, u)$ of problem $\eqref{eq:EVP}$ with $\|u\|_0 = 1$ such that
\begin{align}
|\lambda-\lambda_{h}| & \leq C_{\gamma,\lambda} h^2, \\
\|\mathbf p-\mathbf p_h\|_0 & \leq C_{\gamma,\lambda} h^{3/2}, \label{eq:p_ph L2 estimate EVP divergence-based HDG}\\
\|u-u_h\|_0 & \leq C_{\gamma,\lambda} h.
\end{align}
\end{thrm}
We note that the superconvergence result \eqref{eq:BVP superconvergence result} given for the source problem cannot be generalized in a straightforward way to the eigenvalue problem. But thanks to the very superconvergence result of second order elliptic eigenvalue problems by mixed finite element methods in \cite{MR2891470}, we readily obtain the desired result. To facilitate our analysis, we need to introduce some finite element projectors. Define the HDG finite element projection operator $(G_h, \widetilde{R}_h)$ with $\widetilde{R}_h := (R_h,\hat{R}_h)$ by
\begin{equation*}
G_h \times \widetilde{R}_h: \mathbf Q \times V \to \mathbf Q_h \times \widetilde{V}_h,
\end{equation*}
such that for any $(\mathbf p, u) \in \mathbf Q \times V$
\begin{equation}
\begin{cases}
a_h(G_h(\mathbf p, u), \mathbf q_h)+b_h(\mathbf q_h, \widetilde{R}_h(\mathbf p,u))=a_h(\mathbf p, \mathbf q_h)+b_h(\mathbf q_h, u) &\forall\, \mathbf q_h \in \mathbf Q_h, \\
b_h(G_h(\mathbf p, u), \widetilde{v}_h)+c_h(\widetilde{R}_h(\mathbf p,u), \widetilde{v}_h) = b_h(\mathbf p, \widetilde{v}_h)+c_h(u, \widetilde{v}_h) & \forall\, \widetilde{v}_h \in \widetilde{V}_h.
\end{cases}
\end{equation}
Notice that $G_h(\mathbf p, u)$, $R_h(\mathbf p,u)$, $\hat{R}_h(\mathbf p,u)$ are the components of the projection in $\mathbf Q_h$, $V_h$ and $\hat{V}_h$, respectively.

According to the definition of $R_h$ and Theorem \ref{thm:a priori estimate BVP divergence-based HDG}, we may arrive at the following theorem immediately.
\begin{thrm} \label{thm:R_h(bm p, u) - Q_h0 u}
Let $(\lambda, \mathbf p, u)$ be an eigenpair of problem \eqref{eq:EVP}. Then it holds that 
\begin{equation}
\|R_h(\mathbf p, u) - Q_h^0 u\| \le C_{\gamma,4} h^2.
\end{equation}
\end{thrm}

The following theorem reveals the superconvergence between discrete eigenfunctions and the finite element projection of continuous eigenfunctions. The proof can be obtained by repeating step by step the procedure shown in \cite{MR2891470}. Although the requirement of conforming finite element spaces of mixed form in \cite{MR2891470} is no longer satisfied, the conclusion still holds thanks to the consistency property of the stabilized HDG methods \eqref{eq:consistency}. For brevity, we omit the proof of the following theorem.

\begin{thrm}\label{thm:R_h(bm p, u) - uh}
Let $(\lambda, \mathbf p, u)$ with $\|u\|_0 = 1$ be an eigenpair of problem \eqref{eq:EVP} and $(\lambda_h, \mathbf p_h, \widetilde{u}_h)$ with $\|u_h\|_0 = 1$ be an eigenpair of discrete problem \eqref{eq:EVP HDG discretization}. Then it holds that 
\begin{equation}\label{eq:R_h(bm p, u) - u_h L2 estimate}
\|R_h(\mathbf p, u) - u_h\| \le C_{\gamma} h^2.
\end{equation}
\end{thrm}

Based on Theorem \ref{thm:R_h(bm p, u) - Q_h0 u}, Theorem \ref{thm:R_h(bm p, u) - uh}, Lemma \ref{lem:edge estimate lem} (see the remark below) and \eqref{eq:Qe^1 L2 estimate}, we obtain the following corollary.
\begin{crllr}\label{cor:Qh0 u_uh superconvergence result EVP}
Let $(\lambda, \mathbf p, u)$ with $\|u\|_0 = 1$ be an eigenpair of problem \eqref{eq:EVP} and $(\lambda_h, \mathbf p_h, \widetilde{u}_h)$ with $\|u_h\|_0 = 1$ be an eigenpair of discrete problem \eqref{eq:EVP HDG discretization}. The following superconvergence results hold when $h$ is sufficiently small
\begin{align}
\|Q_h^0 u - u_h\|_0 &\le C_{\gamma} h^2, \\
\|Q_{h,e}^1 u-\hat{u}_h\|_{\mathcal{E}_h^i} &\le C_{\gamma} h^{3/2}, \\
\|u-\hat{u}_h\|_{\mathcal{E}_h^i} &\le C_{\gamma} h^{3/2}.
\end{align}
\end{crllr}

Now, let us state the upper bound properties of eigenvalues for the divergence-based HDG method, which is the main purpose in this section.

\begin{thrm}\label{thm:Upper Bounds for Eigenvalues divergence-based}
For any given penalty parameter $\gamma \ge 0$, it holds that $\lambda \le \lambda_h$ for sufficiently small $h$.
\end{thrm}
\begin{proof}
It follows from \eqref{eq:lam_lamh equality}, the definition of $A_h(\cdot,\cdot)$, the property of projectors $Q_h^0$, $Q_{h,e}^1$ and the continuity of the normal component of $\mathbf p$ at the interelement boundaries that
\begin{align}
\lambda-\lambda_h &= \|\mathbf p-\mathbf p_h\|_{0}^2 - 2(u-u_h,\operatorname{div}(\mathbf p-\mathbf p_h))_{\mathcal{T}_h} + 2\langle u-\hat{u}_h,(\mathbf p-\mathbf p_h)\cdot \mathbf n\rangle_{\partial \mathcal{T}_h} \notag \\
&\quad - \gamma h \|u_h-\hat{u}_h\|_{\partial \mathcal{T}_h}^2 + \lambda \|u-u_h\|_0^2 \notag \\
&= \|\mathbf p-\mathbf p_h\|_{0}^2 - \lambda \|u-u_h\|_0^2 - 2(u-u_h, \lambda u_h-\operatorname{div}\mathbf p_h)_{\mathcal{T}_h} + 2\langle u-\hat{u}_h,(\mathbf p-\mathbf p_h)\cdot \mathbf n\rangle_{\partial \mathcal{T}_h} \notag \\
&\quad - \gamma h \|u_h-\hat{u}_h\|_{\partial \mathcal{T}_h}^2 \notag \\
&= \|\mathbf p-\mathbf p_h\|_{0}^2 - \lambda \|u-u_h\|_0^2 - 2\left(Q_h^0 u-u_h, \lambda u_h-\operatorname{div}\mathbf p_h\right)_{\mathcal{T}_h} - 2\left\langle Q_{h,e}^1 u-\hat{u}_h,[\mathbf p_h]\right\rangle_{\mathcal{E}_h^i} \notag \\
&\quad - \gamma h \|u_h-\hat{u}_h\|_{\partial \mathcal{T}_h}^2. \label{eq:lam_lamh equality for divergence-based HDG}
\end{align}
According to the convergence of eigenfunctions \eqref{eq:p_ph L2 estimate EVP divergence-based HDG}, we are well aware that the term $\|\mathbf p-\mathbf p_h\|_{0}^2$ in the above equalities is a high order term. Due to $- \lambda \|u-u_h\|_0^2 <0$ and $- \gamma h \|u_h-\hat{u}_h\|_{\partial \mathcal{T}_h}^2 <0$, it is sufficient to focus on the estimation of the third and fourth terms. By Corollary \ref{cor:Qh0 u_uh superconvergence result EVP} and a similar deduction of \eqref{eq:wh_hat_wh L2e estimate}, we get
\begin{equation}\label{eq:uh_hat_uh L2e estimate}
\|u_h-\hat{u}_h\|_{\partial \mathcal{T}_h} \le C_\gamma h^{\frac{1}{2}}.
\end{equation}
Choosing $\hat{v}_h = 0$ in the second equation of \eqref{eq:EVP HDG discretization} leads to
\begin{equation*}
(\lambda_h u_h - \operatorname{div} \mathbf p_h, v_h)_{\mathcal{T}_h} = \gamma\langle h_K(u_h-\hat{u}_h), v_h\rangle_{\partial \mathcal{T}_h} \le C\gamma \sum_{K \in \mathcal{T}_h}h_K^{\frac{1}{2}}\|u_h - \hat{u}_h\|_{0,\partial K} \|v_h\|_{0,K} \quad \forall\, v_h \in V_h,
\end{equation*}
which yields
\begin{equation}\label{eq:lamh uh_div ph}
\|\lambda_h u_h - \operatorname{div} \mathbf p_h\|_0 \le C\gamma h^{\frac{1}{2}}\|u_h-\hat{u}_h\|_{\partial \mathcal{T}_h} \le C_\gamma h.
\end{equation}
Similarly, choosing $v_h = 0$ in the second equation of \eqref{eq:EVP HDG discretization} leads to
\begin{align*}
\langle [\mathbf p_h], \hat{v}_h\rangle_{\mathcal{E}_h^i} &= \langle \hat{v}_h, \mathbf p_h \cdot \mathbf n\rangle_{\partial \mathcal{T}_h} = -\gamma\langle h_K(u_h-\hat{u}_h), \hat{v}_h\rangle_{\partial \mathcal{T}_h} \le \gamma h\sum_{K \in \mathcal{T}_h} \|u_h - \hat{u}_h \|_{0, \partial K}\|\hat{v}_h\|_{0,\partial K} \\
&\le C\gamma h \|u_h - \hat{u}_h \|_{\partial \mathcal{T}_h}\|\hat{v}_h\|_{\mathcal{E}_h^i}.
\end{align*}
By \eqref{eq:uh_hat_uh L2e estimate} and the inequalities above, we obtain
\begin{equation}\label{eq:[ph]}
\left\|[\mathbf p_h]\right\|_{\mathcal{E}_h^i} \le C\gamma h \|u_h - \hat{u}_h \|_{\partial \mathcal{T}_h} \le C_\gamma h^{\frac{3}{2}}.
\end{equation}

Finally combining Lemma \ref{lem:The Saturation Condition}, Theorem \ref{thm:Convergence of eigenvalues and eigenfunctions divergence-based HDG}, Corollary \ref{cor:Qh0 u_uh superconvergence result EVP}, \eqref{eq:lam_lamh equality for divergence-based HDG}, \eqref{eq:lamh uh_div ph}, \eqref{eq:[ph]}, we get
\begin{align*}
\lambda-\lambda_h &\le \|\mathbf p-\mathbf p_h\|_{0}^2 - \lambda \|u-u_h\|_0^2 - 2\left(Q_h^0 u-u_h, \lambda u_h-\operatorname{div}\mathbf p_h\right)_{\mathcal{T}_h} - 2\left\langle Q_{h,e}^1 u-\hat{u}_h,[\mathbf p_h]\right\rangle_{\mathcal{E}_h^i} \\
&\le C_\gamma h^3 -C h^2 + 2\left\|Q_h^0 u-u_h\right\|_0 \left\|\lambda u_h-\operatorname{div}\mathbf p_h\right\|_0 + 2\left\| Q_{h,e}^1 u-\hat{u}_h\right\|_{\mathcal{E}_h^i} \left\|[\mathbf p_h]\right\|_{\mathcal{E}_h^i} \\
&\le C_\gamma h^3 -C h^2 + 2\left\|Q_h^0 u-u_h\right\|_0 \left\|\lambda u_h - \lambda_h u_h + \lambda_h u_h-\operatorname{div}\mathbf p_h\right\|_0 \\
&\le C_\gamma h^3 -C h^2 + 2\left\|Q_h^0 u-u_h\right\|_0 \left(|\lambda - \lambda_h| \|u_h\|_0 + \left\|\lambda_h u_h-\operatorname{div}\mathbf p_h\right\|_0 \right) \\
&\le C_\gamma h^3 -C h^2,
\end{align*}
which indicates that $\lambda-\lambda_h$ is negative when $h$ is sufficiently small.
\end{proof}

One can easily derive the following corollary based on the aforementioned theorem, as the divergence-based HDG method degenerates into the BDM-H method for $\gamma = 0$ and notice that the original BDM method is a BDM-H method; see \cite{MR2485455}.
\begin{crllr}\label{cor:Upper Bounds for Eigenvalues BDM}
Suppose that $\lambda_h$ is the $k$th discrete eigenvalue of the Laplace operator computed by the BDM element (see \cite{MR799685}), $\lambda$ is the $k$th eigenvalue of \eqref{eq:EVP}, and the corresponding eigenfunctions $\mathbf p \in \mathbf H^2(\Omega)$, $\operatorname{div} \mathbf p \in H^1(\Omega)$, $u \in H^1(\Omega)$. Then it holds that $\lambda \le \lambda_h$ for sufficiently small $h$.
\end{crllr}

\section{A High Accuracy Algorithm} \label{sec:A High Accuracy Algorithm}

It is apparent that the gradient-based HDG method yields both upper and lower bounds of eigenvalues, whereas the divergence-based HDG method can only compute upper bounds. Therefore, we could take advantage of the first type of HDG methods, in conjunction with a post-processing technique, to obtain a higher accuracy approximation. Assume that $\{\lambda^{L}_h\}$ and $\{\lambda^{U}_h\}$ computed by the HDG method with the penalty parameter $\gamma_1$ and $\gamma_2$, respectively, satisfy
\begin{equation}
\lambda^{L}_{h_1} \le \lambda^{L}_{h_2} \cdots \le \lambda^{L}_{h_n} \le \lambda \le \lambda^{U}_{h_n} \le \cdots \le \lambda^{U}_{h_2} \le \lambda^{U}_{h_1},
\end{equation}
where $h_1 \ge h_2 \ge \cdots \searrow 0$. Therefore, we have
\begin{align*}
	\lambda^{L}_{h_n} = \lambda - C_1 h_n^2, \qquad\lambda^{U}_{h_n} = \lambda + C_2 h_n^2,
\end{align*}
where $C_i > 0\ (i = 1,2)$. Theoretically speaking, we may accelerate the convergence via
\begin{equation*}
\widehat\lambda_h = \rho \lambda^{U}_h + (1 - \rho) \lambda^{L}_h,
\end{equation*}
with $\rho = \frac{C_1}{C_1+C_2} > 0$. Since the constants $C_1$ and $C_2$ are usually unavailable in practical computations, we may compute $\rho$ asymptotically with
\begin{equation}
	\rho_{h_{n+1}} = \frac{\lambda^{L}_{h_n} - \lambda^{L}_{h_{n+1}}}{\lambda^{U}_{h_{n+1}} - \lambda^{U}_{h_n} + \lambda^{L}_{h_n} - \lambda^{L}_{h_{n+1}}}.
\end{equation}
Thus the convex combination of two discrete eigenvalues sequences may be written as
\begin{equation}
	\widehat\lambda_{h_{n+1}} = \rho_{h_{n+1}} \lambda^{U}_{h_{n+1}} + (1 - \rho_{h_{n+1}}) \lambda^{L}_{h_{n+1}}.
\end{equation}
The implementation of such post-processing technique has been found to be a highly effective approach in eliminating dominant errors and achieving significantly improved accuracy. The subsequent section presents numerical findings that validate the aforementioned high accuracy approximation.

\section{Numerical Experiments}
\label{sec:Numerical Experiments}
In this section, some numerical examples are presented to verify that
\begin{enumerate}[(1)]
\item The gradient-based HDG method can compute both upper and lower bounds of exact eigenvalues through adjusting penalty parameter carefully.

\item Both the divergence-based HDG method and the BDM mixed finite element method can provide upper bounds of eigenvalues.

\item High accuracy post-processing algorithm can be applied to improve convergence order and accuracy using the gradient-based HDG finite element.
\end{enumerate}

If we know the true continuous eigenvalue $\lambda$, then the convergence rate is computed by
\begin{equation*}
\text{Ratio} =\log_2\left(\frac{\left|\lambda-\lambda_{2 h}\right|}{\left|\lambda-\lambda_h\right|}\right).
\end{equation*}
 Otherwise,
\begin{equation*}
\text{Ratio} =\log_2\left(\frac{\left|\lambda_{4 h}-\lambda_{2 h}\right|}{\left|\lambda_{2 h}-\lambda_h\right|}\right).
\end{equation*}
The initial mesh is obtained by splitting the domain $\Omega$  into two right triangles by its positively sloped diagonal, and each triangle is subdivided into four smaller triangles uniformly to get a finer triangulation. All computations are performed on a reduced nonlinear eigenvalue problem using the technique of static condensation.

Subsection \ref{subs:Laplacian Eigenvalue Problem} computes eigenvalues of the Laplace operator on a square domain with elements of the lowest degree. Subsection \ref{subs:High Order Element} computes eigenvalues with higher order polynomials. Subsection \ref{subs:High Discontinuous Coefficients Eigenvalue Problem} deals with eigenvalues of a second order elliptic operator with high discontinuous coefficients.

\subsection{Laplacian Eigenvalue Problem}
\label{subs:Laplacian Eigenvalue Problem}
We consider the Laplacian eigenvalue problem \eqref{eq:EVP} on domain $\Omega = (0,\pi)^2$. Its solutions are
\begin{equation*}
	\lambda = m^2 + n^2, \quad u = \sin(mx)\sin(ny),
\end{equation*}
where $m$, $n$ are arbitrary integers.

\begin{table}
\centering
\setlength{\abovecaptionskip}{0cm}
\setlength{\belowcaptionskip}{0.2cm}
\caption{The first ten discrete eigenvalues by the gradient-based HDG method with $\gamma=1$, $10$}
\label{tab:The first ten eig gamma=1 and 10 HDG}
\begin{tabular}{crrrrrr}
\toprule
& \multicolumn{3}{c}{$\gamma=1$} & \multicolumn{3}{c}{$\gamma=10$} \\
\cmidrule(lr){2-4}\cmidrule(lr){5-7}
Exact & $h=\frac{\pi}{2^6}$ & $h=\frac{\pi}{2^7}$ & $h=\frac{\pi}{2^8}$ & $h=\frac{\pi}{2^6}$ & $h=\frac{\pi}{2^7}$ & $h=\frac{\pi}{2^8}$ \\ 
\midrule
2  & 1.997896  & 1.999474  & 1.999868  & 2.000174  & 2.000043  & 2.000011  \\
5  & 4.986194  & 4.996541  & 4.999135  & 5.000437  & 5.000109  & 5.000027  \\
5  & 4.986310  & 4.996570  & 4.999142  & 5.001390  & 5.000348  & 5.000087  \\
8  & 7.966440  & 7.991583  & 7.997894  & 8.002777  & 8.000695  & 8.000174  \\
10 & 9.943253  & 9.985751  & 9.996434  & 10.003115 & 10.000780 & 10.000195 \\
10 & 9.943253  & 9.985751  & 9.996434  & 10.003119 & 10.000780 & 10.000195 \\
13 & 12.909430 & 12.977237 & 12.994302 & 13.002837 & 13.000712 & 13.000178 \\
13 & 12.910404 & 12.977482 & 12.994363 & 13.010865 & 13.002720 & 13.000680 \\
17 & 16.833928 & 16.958167 & 16.989522 & 17.007977 & 17.001996 & 17.000499 \\
17 & 16.833990 & 16.958183 & 16.989526 & 17.008495 & 17.002126 & 17.000532 \\ 
\midrule
\makecell{ DOF \\(condensed) } & 24320 & 97792 & 392192 & 24320 & 97792 & 392192 \\
\bottomrule
\end{tabular}
\end{table}

Table \ref{tab:The first ten eig gamma=1 and 10 HDG} presents the first ten eigenvalues computed by the gradient-based HDG method of the lowest degree. The first ten eigenvalues, obtained with the stabilization parameter $\gamma$ (taking 1 and 10), are compared with the exact ones on the last three successively refined Cartesian grids. In the last row, DOF denotes the total degrees of freedom of the reduced system. For $\gamma = 1$, the computed eigenvalues approximate the exact ones from below, while for $\gamma = 10$, the approximation is from above. These results verify Theorem \ref{thm:Lower Bounds for Eigenvalues gradient-based} and Theorem \ref{thm:Upper Bounds for Eigenvalues gradient-based}. Although $\gamma$ is theoretically assumed to be sufficiently large (or small), we find that it suffices to ensure that $\gamma \geq 10$ (or $\gamma \leq 1$). The results obtained for $\gamma > 10$ (or $\gamma < 1$) are completely analogous and are not reported, for the sake of brevity.

\begin{figure}
 \centering
 \includegraphics[width=8cm]{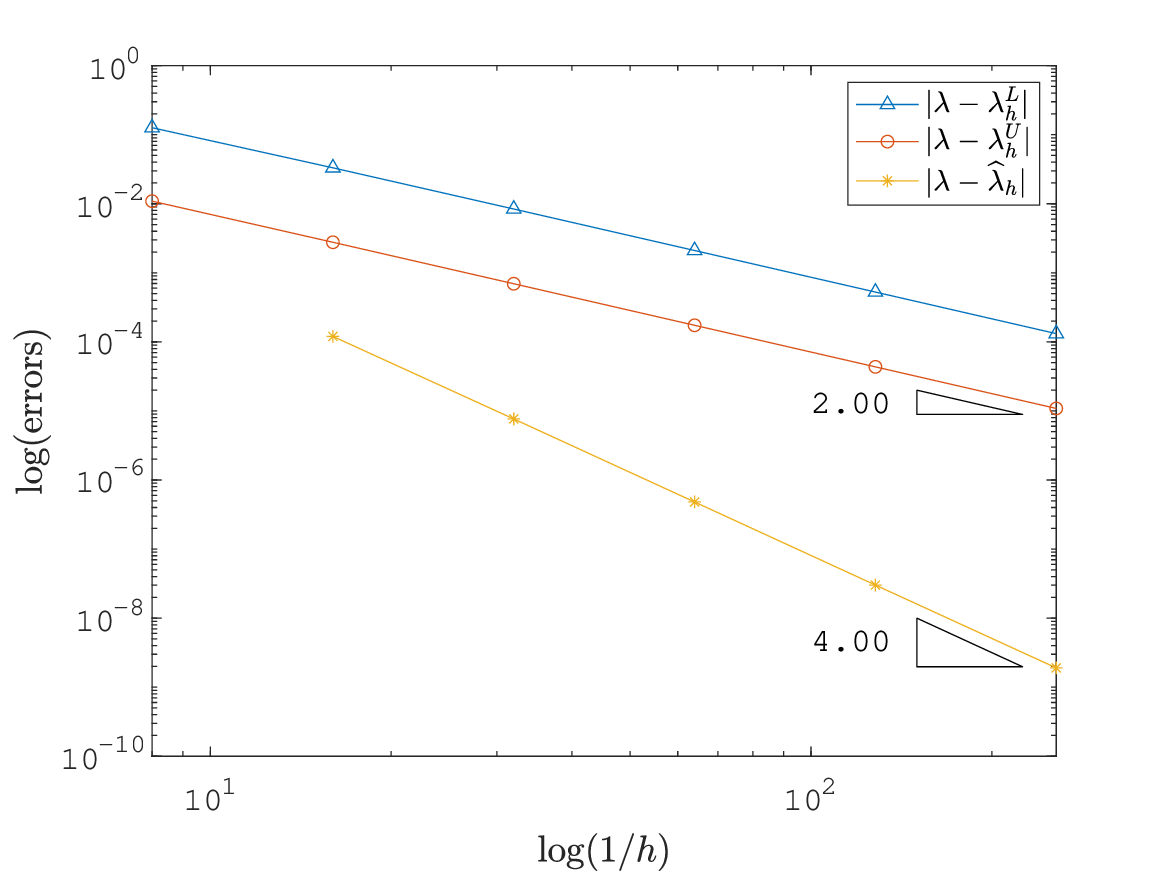}
 \caption{The errors of the first eigenvalue $\lambda_h^L, \lambda_h^U$ and $\widehat\lambda_h$}
 \label{fig:err_1st_eigenvalue}
\end{figure}

We now turn to the application of post-processing techniques. Table \ref{tab:1st convex combination gamma=1 and 10 HDG} -- \ref{tab:10th convex combination gamma=1 and 10 HDG} present the discrete eigenvalues $\lambda_{1,h}$, $\lambda_{2,h}$ and $\lambda_{10,h}$ computed by the gradient-based HDG method, as well as the approximate eigenvalues $\widehat{\lambda}_{h}$ computed by the high accuracy algorithm. We observe that the algorithm improves the convergence rate of discrete eigenvalues from $O(h^2)$ to $O(h^4)$ and a better accuracy is achieved. Figure \ref{fig:err_1st_eigenvalue} (log-log scale) plots the errors of the first eigenvalue $\lambda_h^L$, $\lambda_h^U$ and $\widehat\lambda_h$. We clearly observe quadratic convergence of discrete eigenvalues by the HDG method (as predicted by Theorem \ref{thm:Convergence of eigenvalues and eigenfunctions gradient-based}), and quartic convergence of approximate eigenvalues by the post-processing technique.

\begin{table}
\centering
\setlength{\abovecaptionskip}{0cm}
\setlength{\belowcaptionskip}{0.2cm}
\caption{The first eigenvalue by the gradient-based HDG method with $\gamma_1=1$ and $\gamma_2=10$}
\label{tab:1st convex combination gamma=1 and 10 HDG}
\begin{tabular}{ccccccc}
\toprule
$h/\pi$&   $\lambda^L_h$ & Ratio  & $\lambda^U_h$ & Ratio &    $\widehat\lambda_h$      &     Ratio   \\ \midrule
$2^{-3}$ & 1.8738170118 & -      & 2.0109072750 & -      & -            & -      \\
$2^{-4}$ & 1.9668630165 & 1.9290 & 2.0027681706 & 1.9783 & 1.9998800401 & -      \\
$2^{-5}$ & 1.9916100480 & 1.9817 & 2.0006946768 & 1.9945 & 1.9999923445 & 3.9699 \\
$2^{-6}$ & 1.9978958017 & 1.9954 & 2.0001738349 & 1.9986 & 1.9999995194 & 3.9936 \\
$2^{-7}$ & 1.9994735294 & 1.9988 & 2.0000434691 & 1.9997 & 1.9999999699 & 3.9992 \\
$2^{-8}$ & 1.9998683560 & 1.9997 & 2.0000108679 & 1.9999 & 1.9999999981 & 4.0024 \\ 
\bottomrule
\end{tabular}
\end{table}

\begin{table}
\centering
\setlength{\abovecaptionskip}{0cm}
\setlength{\belowcaptionskip}{0.2cm}
\caption{The second eigenvalue by the gradient-based HDG method with $\gamma_1=1$ and $\gamma_2=10$}
\label{tab:2nd convex combination gamma=1 and 10 HDG}
\begin{tabular}{ccccccc}
\toprule
$h/\pi$&   $\lambda^L_h$ & Ratio  & $\lambda^U_h$ & Ratio &    $\widehat\lambda_h$      &     Ratio   \\ \midrule
$2^{-3}$ & 4.2497864709 & -      & 5.0262947520 & -      & -            & -      \\
$2^{-4}$ & 4.7880598769 & 1.8236 & 5.0068848547 & 1.9333 & 4.9992687573 & -      \\
$2^{-5}$ & 4.9452388991 & 1.9524 & 5.0017421129 & 1.9826 & 4.9999519559 & 3.9279 \\
$2^{-6}$ & 4.9861941190 & 1.9879 & 5.0004368644 & 1.9956 & 4.9999969658 & 3.9849 \\
$2^{-7}$ & 4.9965412301 & 1.9970 & 5.0001093002 & 1.9989 & 4.9999998101 & 3.9979 \\
$2^{-8}$ & 4.9991348502 & 1.9992 & 5.0000273303 & 1.9997 & 4.9999999882 & 4.0024 \\
\bottomrule
\end{tabular}
\end{table}

\begin{table}
\centering
\setlength{\abovecaptionskip}{0cm}
\setlength{\belowcaptionskip}{0.2cm}
\caption{The 10th eigenvalue by the gradient-based HDG method with $\gamma_1=1$ and $\gamma_2=10$}
\label{tab:10th convex combination gamma=1 and 10 HDG}
\begin{tabular}{ccccccc}
\toprule
$h/\pi$&   $\lambda^L_h$ & Ratio  & $\lambda^U_h$ & Ratio &    $\widehat\lambda_h$      &     Ratio   \\ \midrule
$2^{-3}$ & 10.4967125736 & -      & 17.4978310045 & -      & -             & -      \\
$2^{-4}$ & 14.6915529845 & 1.4942 & 17.1331677984 & 1.9024 & 16.9378906014 & -      \\
$2^{-5}$ & 16.3553674341 & 1.8404 & 17.0338387825 & 1.9765 & 16.9956162029 & 3.8246 \\
$2^{-6}$ & 16.8339903595 & 1.9572 & 17.0084946209 & 1.9941 & 16.9997189204 & 3.9631 \\
$2^{-7}$ & 16.9581829950 & 1.9891 & 17.0021258589 & 1.9985 & 16.9999823338 & 3.9919 \\
$2^{-8}$ & 16.9895259035 & 1.9973 & 17.0005316030 & 1.9996 & 16.9999988947 & 3.9985 \\
\bottomrule
\end{tabular}
\end{table}

In Table \ref{tab:The first ten eig divergence-based gamma=0 and 50 HDG}, we report the discrete eigenvalues by the divergence-based HDG method with $\gamma=0$ and $50$. For brevity, only the values computed on the last three successively refined Cartesian grids are reported. A simple calculation shows that both the BDM method ($\gamma=0$) and the divergence-based HDG method ($\gamma=50$) exhibit expected quadratic convergence. We observe that the discrete eigenvalues obtained by both methods approach the exact ones from above, which is consistent with Theorem \ref{thm:Upper Bounds for Eigenvalues divergence-based} and Corollary \ref{cor:Upper Bounds for Eigenvalues BDM}.

\begin{table}
\centering
\setlength{\abovecaptionskip}{0cm}
\setlength{\belowcaptionskip}{0.2cm}
\caption{The first ten discrete eigenvalues by the divergence-based HDG method with $\gamma=0$, $50$}
\label{tab:The first ten eig divergence-based gamma=0 and 50 HDG}
\begin{tabular}{crrrrrr}
\toprule
& \multicolumn{3}{c}{$\gamma=0$} & \multicolumn{3}{c}{$\gamma=50$} \\
\cmidrule(lr){2-4}\cmidrule(lr){5-7}
Exact & $h=\frac{\pi}{2^6}$ & $h=\frac{\pi}{2^7}$ & $h=\frac{\pi}{2^8}$ & $h=\frac{\pi}{2^6}$ & $h=\frac{\pi}{2^7}$ & $h=\frac{\pi}{2^8}$ \\
\midrule
2  & 2.000535  & 2.000134  & 2.000033  & 2.249280  & 2.064117  & 2.016151  \\
5  & 5.002864  & 5.000716  & 5.000179  & 5.534974  & 5.137611  & 5.034665  \\
5  & 5.003829  & 5.000957  & 5.000239  & 5.715131  & 5.183956  & 5.046339  \\
8  & 8.008563  & 8.002142  & 8.000535  & 8.977158  & 8.256151  & 8.064880  \\
10 & 10.013386 & 10.003347 & 10.000837 & 11.255084 & 10.323120 & 10.081415 \\
10 & 10.013390 & 10.003347 & 10.000837 & 11.283347 & 10.325182 & 10.081549 \\
13 & 13.018551 & 13.004640 & 13.001160 & 14.328465 & 13.344578 & 13.087027 \\
13 & 13.026679 & 13.006671 & 13.001668 & 14.914303 & 13.495922 & 13.125160 \\
17 & 17.038433 & 17.009607 & 17.002402 & 19.151057 & 17.550775 & 17.138538 \\
17 & 17.038958 & 17.009738 & 17.002434 & 19.182997 & 17.558339 & 17.140429 \\ 
\midrule
\makecell{ DOF \\(condensed) } & 24320 & 97792 & 392192 & 24320 & 97792 & 392192 \\
\bottomrule
\end{tabular}
\end{table}

\subsection{High Order Element}
\label{subs:High Order Element}
\begin{figure}
 \centering
 \includegraphics[ width=8cm]{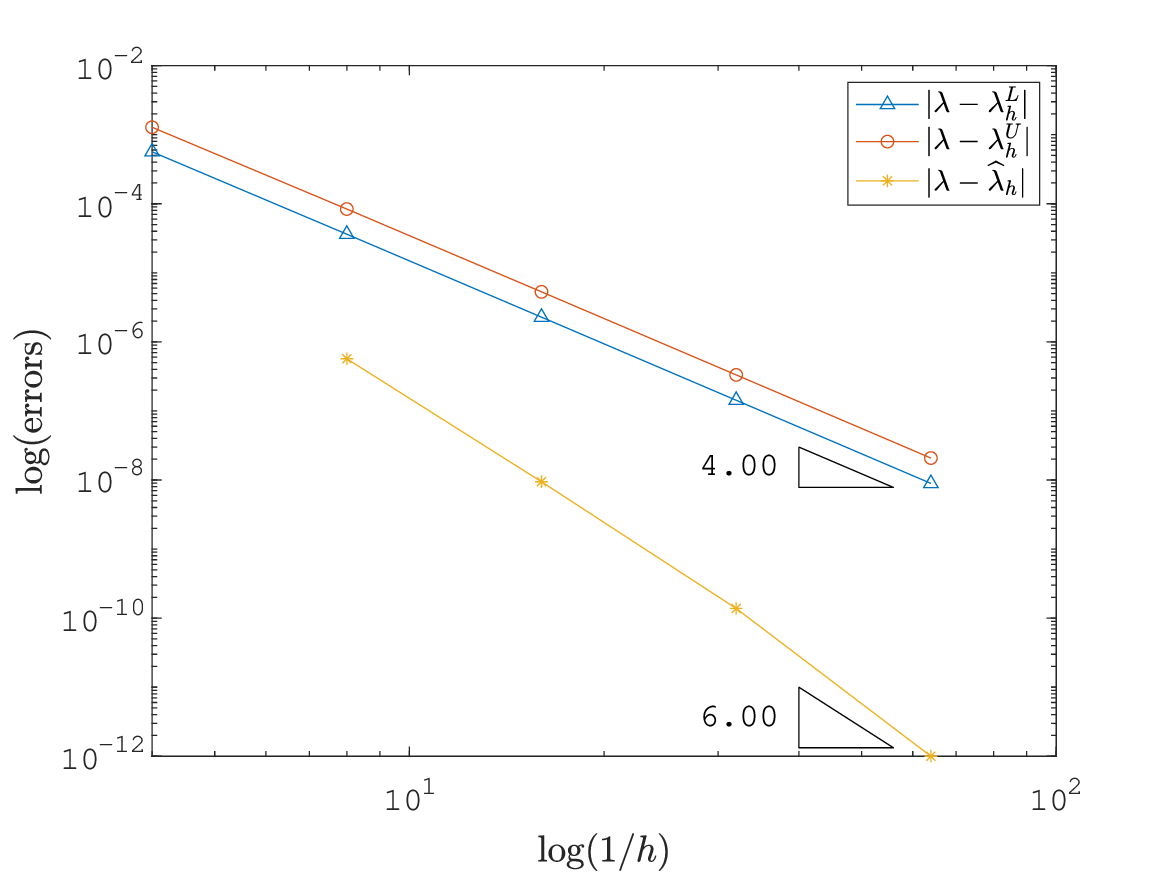}
 \caption{The errors of the first eigenvalue $\lambda_h^L, \lambda_h^U$ and $\widehat\lambda_h$, $k = 2$}
 \label{fig:err_1st_eigenvalue_P2}
\end{figure}

In this subsection, we apply the gradient-based HDG method with higher order polynomials to the mixed Laplacian eigenvalue problem \eqref{eq:EVP} on $\Omega = (0,\pi)^2$. The corresponding local spaces are defined as:
\begin{equation*}
V(K) = \mathcal{P}_2(K), \quad \mathbf Q(K) = \bm{\mathcal{P}}_1(K), \quad \hat{V}(e) = \mathcal{P}_2(e),
\end{equation*}
and stabilization parameters $\gamma_1=8$, $\gamma_2=15$ are employed. Table \ref{tab:1st convex combination gamma=8 and 15 HDG P2} -- \ref{tab:10th convex combination gamma=8 and 15 HDG P2} show that the HDG method produces lower ($\lambda^L_h$) and upper ($\lambda^U_h$) bounds of eigenvalues with $\gamma_1$, $\gamma_2$, respectively, and the quartic convergence is obtained. In addition, the convergence rate of discrete eigenvalues is also improved by the post-processing algorithm. Figure \ref{fig:err_1st_eigenvalue_P2} (log-log scale) plots the errors of the first eigenvalue $\lambda_h^L$, $\lambda_h^U$ and the corresponding approximate eigenvalue $\widehat\lambda_h$. We clearly observe quartic convergence of discrete eigenvalues by the HDG method, and higher convergence of approximate eigenvalues by the post-processing technique. 

Concerning the divergence-based HDG method with higher order polynomials, the same tests are repeated with stabilization parameter $\gamma \ge 0$, and quartic convergence of discrete eigenvalues is obtained. We also observe upper bound properties of computed eigenvalues. For brevity, these results are not reported.

\begin{table}
\centering
\setlength{\abovecaptionskip}{0cm}
\setlength{\belowcaptionskip}{0.2cm}
\caption{The first eigenvalue by the gradient-based HDG method with $\gamma_1=8$ and $\gamma_2=15$, $k=2$}
\label{tab:1st convex combination gamma=8 and 15 HDG P2}
\begin{tabular}{ccccccc}
\toprule
$h/\pi$&   $\lambda^L_h$ & Ratio  & $\lambda^U_h$ & Ratio &    $\widehat\lambda_h$      &     Ratio   \\ \midrule
$2^{-2}$ & 1.999434654974 & -      & 2.001280706765 & -      & -              & -      \\
$2^{-3}$ & 1.999963751394 & 3.9631 & 2.000083854158 & 3.9329 & 2.000000569366 & -      \\
$2^{-4}$ & 1.999997719682 & 3.9906 & 2.000005304363 & 3.9826 & 2.000000009435 & 5.9151 \\
$2^{-5}$ & 1.999999857237 & 3.9975 & 2.000000332521 & 3.9957 & 2.000000000138 & 6.0908 \\
$2^{-6}$ & 1.999999991073 & 3.9993 & 2.000000020796 & 3.9991 & 2.000000000001 & 6.8743 \\
\bottomrule
\end{tabular}
\end{table}

\begin{table}
\centering
\setlength{\abovecaptionskip}{0cm}
\setlength{\belowcaptionskip}{0.2cm}
\caption{The second eigenvalue by the gradient-based HDG method with $\gamma_1=8$ and $\gamma_2=15$, $k=2$}
\label{tab:2nd convex combination gamma=8 and 15 HDG P2}
\begin{tabular}{ccccccc}
\toprule
$h/\pi$&   $\lambda^L_h$ & Ratio  & $\lambda^U_h$ & Ratio &    $\widehat\lambda_h$      &     Ratio   \\ \midrule
$2^{-2}$ & 4.986201967690 & -      & 5.007163399257 & -      & -              & -      \\
$2^{-3}$ & 4.999072649167 & 3.8952 & 5.000515943194 & 3.7954 & 5.000024388423 & -      \\
$2^{-4}$ & 4.999940799044 & 3.9694 & 5.000033269657 & 3.9549 & 5.000000228241 & 6.7395 \\
$2^{-5}$ & 4.999996276816 & 3.9910 & 5.000002093847 & 3.9900 & 5.000000001023 & 7.8010 \\
$2^{-6}$ & 4.999999766849 & 3.9972 & 5.000000131046 & 3.9980 & 4.999999999950 & 4.3487 \\
\bottomrule
\end{tabular}
\end{table}


\begin{table}
\centering
\setlength{\abovecaptionskip}{0cm}
\setlength{\belowcaptionskip}{0.2cm}
\caption{The 10th eigenvalue by the gradient-based HDG method with $\gamma_1=8$ and $\gamma_2=15$, $k=2$}
\label{tab:10th convex combination gamma=8 and 15 HDG P2}
\begin{tabular}{ccccccc}
\toprule
$h/\pi$&   $\lambda^L_h$ & Ratio  & $\lambda^U_h$ & Ratio &    $\widehat\lambda_h$      &     Ratio   \\ \midrule
$2^{-2}$ & 16.624493937059 & -      & 17.074104164135 & -      & -               & -      \\
$2^{-3}$ & 16.957653463663 & 3.1485 & 17.023440388063 & 1.6606 & 17.014756668281 & -      \\
$2^{-4}$ & 16.997145473436 & 3.8909 & 17.001540031675 & 3.9280 & 16.999972371107 & 9.0610 \\
$2^{-5}$ & 16.999817825033 & 3.9699 & 17.000097276530 & 3.9847 & 16.999999300922 & 5.3046 \\
$2^{-6}$ & 16.999988546293 & 3.9914 & 17.000006091882 & 3.9971 & 16.999999983245 & 5.3827 \\
\bottomrule
\end{tabular}
\end{table}

\subsection{High Discontinuous Coefficients Eigenvalue Problem}
\label{subs:High Discontinuous Coefficients Eigenvalue Problem}
Consider a discontinuous coefficients case, which reads
\begin{equation}\label{eq:discontinuous coefficients case}
\left\{\begin{array}{ll}
-\nabla \cdot (\alpha(x) \nabla u)=\lambda u &\text { in } \Omega, \\
u=0 &\text { on } \partial \Omega,
\end{array}\right.
\end{equation}
where $\Omega = (0,1)^2$ and $\alpha (x)$ is a piecewise constant function
\begin{equation}
	\alpha (x) = \begin{cases}
		1    &\text{in } \Omega_1\text{ or }\Omega_4, \\
		10^8 &\text{in } \Omega_2\text{ or }\Omega_3.
	\end{cases}
\end{equation}

\begin{figure}
\centering
\begin{tikzpicture}
\coordinate (A) at (3/4,3/4);
\coordinate (B) at (9/4,3/4);
\coordinate (C) at (3/4,9/4);
\coordinate (D) at (9/4,9/4);
\draw [black](0,0) rectangle (3,3);
\draw [black](3/2,0) -- (3/2,3);
\draw [black](0,3/2) -- (3,3/2);
\draw (A) node{$\Omega_3$};
\draw (B) node{$\Omega_4$};
\draw (C) node{$\Omega_1$};
\draw (D) node{$\Omega_2$};
\end{tikzpicture}
\caption{The partition of domain $\Omega$}
\label{fig:The partition of domain}
\end{figure}
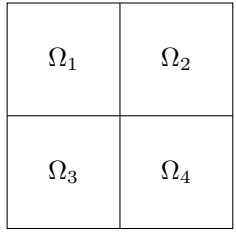

For the gradient-based HDG method of the lowest degree, Table \ref{tab:1st convex combination gamma=5 and 10 HDG high discontinuous coefficients} -- \ref{tab:10th convex combination gamma=5 and 10 HDG high discontinuous coefficients} show that upper and lower bound properties of eigenvalues still hold, i.e., when $\gamma_1 = 5$, the computed eigenvalues approximate the exact ones from below; when $\gamma_2 = 10$, the computed eigenvalues approximate the exact ones from above. The convergence rate of discrete eigenvalues is also improved by the post-processing algorithm.

\begin{table}
\centering
\setlength{\abovecaptionskip}{0cm}
\setlength{\belowcaptionskip}{0.2cm}
\caption{The first eigenvalue by the gradient-based HDG method with $\gamma_1=5$ and $\gamma_2=10$}
\label{tab:1st convex combination gamma=5 and 10 HDG high discontinuous coefficients}
\begin{tabular}{ccccccc}
\toprule
$h$&   $\lambda^L_h$ & Ratio  & $\lambda^U_h$ & Ratio &    $\widehat\lambda_h$      &     Ratio   \\ \midrule
$2^{-3}$ & 76.4586534739 & -      & 80.2255564400 & -      & -             & -      \\
$2^{-4}$ & 78.3449734022 & -      & 79.3457063720 & -      & 79.0273980349 & -      \\
$2^{-5}$ & 78.8065169920 & 2.0310 & 79.0610346330 & 1.6280 & 78.9639393622 & -      \\
$2^{-6}$ & 78.9196565673 & 2.0284 & 78.9836305749 & 1.8788 & 78.9576425756 & 3.3331 \\
$2^{-7}$ & 78.9475943934 & 2.0178 & 78.9636186782 & 1.9516 & 78.9569309163 & 3.1454 \\
$2^{-8}$ & 78.9545312815 & 2.0099 & 78.9585404380 & 1.9785 & 78.9568459526 & 3.0663 \\
$2^{-9}$ & 78.9562592960 & 2.0052 & 78.9572619228 & 1.9899 & 78.9568355596 & 3.0312 \\
\bottomrule
\end{tabular}
\end{table}

\begin{table}
\centering
\setlength{\abovecaptionskip}{0cm}
\setlength{\belowcaptionskip}{0.2cm}
\caption{The second eigenvalue by the gradient-based HDG method with $\gamma_1=5$ and $\gamma_2=10$}
\label{tab:2nd convex combination gamma=5 and 10 HDG high discontinuous coefficients}
\begin{tabular}{ccccccc}
\toprule
$h$&   $\lambda^L_h$ & Ratio  & $\lambda^U_h$ & Ratio &    $\widehat\lambda_h$      &     Ratio   \\ \midrule
$2^{-3}$ & 76.4633983311 & -      & 80.2321767109 & -      & -             & -      \\
$2^{-4}$ & 78.3452842168 & -      & 79.3461330788 & -      & 79.0257507298 & -      \\
$2^{-5}$ & 78.8065373934 & 2.0285 & 79.0610624126 & 1.6361 & 78.9638423954 & -      \\
$2^{-6}$ & 78.9196585275 & 2.0277 & 78.9836330019 & 1.8804 & 78.9576372438 & 3.3186 \\
$2^{-7}$ & 78.9475951954 & 2.0176 & 78.9636195096 & 1.9519 & 78.9569312632 & 3.1358 \\
$2^{-8}$ & 78.9545320111 & 2.0098 & 78.9585411698 & 1.9785 & 78.9568466540 & 3.0607 \\
$2^{-9}$ & 78.9562600210 & 2.0052 & 78.9572626487 & 1.9899 & 78.9568362834 & 3.0283 \\
\bottomrule
\end{tabular}
\end{table}

\begin{table}
\centering
\setlength{\abovecaptionskip}{0cm}
\setlength{\belowcaptionskip}{0.2cm}
\caption{The 10th eigenvalue by the gradient-based HDG method with $\gamma_1=5$ and $\gamma_2=10$}
\label{tab:10th convex combination gamma=5 and 10 HDG high discontinuous coefficients}
\begin{tabular}{ccccccc}
\toprule
$h$&   $\lambda^L_h$ & Ratio  & $\lambda^U_h$ & Ratio &    $\widehat\lambda_h$      &     Ratio   \\ \midrule
$2^{-3}$ & 336.6615417977 & -      & 412.4590860207 & -      & -              & -      \\
$2^{-4}$ & 377.2972199614 & -      & 401.6327304065 & -      & 396.5131330880 & -      \\
$2^{-5}$ & 390.1786028321 & 1.6575 & 396.6668059226 & 1.1244 & 394.8614956904 & -      \\
$2^{-6}$ & 393.6190849212 & 1.9046 & 395.2684807493 & 1.8284 & 394.7918360759 & 4.5674 \\
$2^{-7}$ & 394.4922947590 & 1.9782 & 394.9064559881 & 1.9495 & 394.7850728757 & 3.3645 \\
$2^{-8}$ & 394.7111973951 & 1.9960 & 394.8148636630 & 1.9828 & 394.7842833454 & 3.0986 \\
$2^{-9}$ & 394.7659317713 & 1.9998 & 394.7918568770 & 1.9932 & 394.7841845791 & 2.9989 \\
\bottomrule
\end{tabular}
\end{table}

\section{Conclusions}
\label{sec:Conclusions}
In this paper, we apply two types of HDG methods to solve elliptic eigenvalue problems, and develop corresponding convergence analysis for eigenvalues and eigenfunctions. For the gradient-based HDG method, we prove that if $\gamma$ is sufficiently large (resp. small), discrete eigenvalues by the HDG method can approximate exact eigenvalues from above (resp. below); for the divergence-based HDG method, we prove the discrete eigenvalues approximate the exact ones from above no matter how the penalty parameter changes. Furthermore, upper bound properties of the BDM element are also established. Finally, a very effective post-processing algorithm is designed to obtain higher convergence order. Numerical experiments confirm our theoretical findings.
\bibliographystyle{plain}
\bibliography{HDG}
\end{document}